%% file: main.tex
\documentclass[preprint,12pt]{elsarticle}
\usepackage{bm}
\usepackage{amsmath,amsfonts,amssymb, xcolor}

\input{common.tex}
\usepackage{array, tabularx, multirow, booktabs, caption}
\usepackage{pgffor, pgfplotstable}
\pgfplotsset{compat=1.18} 
\usepackage{enumitem}

\journal{Elsevier}

\newcommand{\bn}{\ensuremath{\mathbf n}}
\newcommand{\bu}{\ensuremath{\mathbf u}}
\newcommand{\vWERPfirst}{vWERP\ensuremath{_\mathrm{1st}}}
\newcommand{\vWERPsecond}{vWERP\ensuremath{_\mathrm{2nd}}}

\begin{document}

\begin{frontmatter}
%\title{Physics Informed Neural Networks State Estimation of the Aortic Blood Flow Including Hematocrit Dependent Rheological Behavior}
\title{Estimation of Hemodynamic Parameters via Physics Informed Neural Networks including Hematocrit Dependent Rheology}
\author[inst1,inst2]{Moises Sierpe}
\affiliation[inst1]{organization={Universidad de Santiago de Chile, Av Libertador Bernardo O'Higgins 3363, Santiago, Chile}}
\ead{moises.sierpe@usach.cl} 
\affiliation[inst2]{Computational Heat and Fluid Flow Lab. Universidad de Santiago de Chile.}
\author[inst1,inst2]{Ernesto Castillo}
\ead{ernesto.castillode@usach.cl}
\author[inst3]{Hernán Mella}
\ead{hernan.mella@pucv.cl}
\author[inst4]{Felipe Galarce$^*$\footnote{$^*$Corresponding author}}
\ead{felipe.galarce@pucv.cl}
\affiliation[inst3]{organization={School of Electrical Engineering. Pontificia Universidad Católica de Valparaíso. Avenida Brasil 2147, Valparaíso. Chile.}}
\affiliation[inst4]{organization={School of Civil Engineering. Pontificia Universidad Católica de Valparaíso. Avenida Brasil 2147, Valparaíso. Chile.}}

\begin{abstract}
Physics-Informed Neural Networks (PINNs) show significant potential for solving inverse problems, especially when observations are limited and sparse, provided that the relevant physical equations are known. We use PINNs to estimate smooth velocity and pressure fields from synthetic 4D flow Magnetic Resonance Imaging (MRI) data. We analyze five non-Newtonian dynamic 3D blood flow cases within a realistic aortic model, covering a range of hematocrit values from anemic to polycythemic conditions. To enhance state estimation results, we consider various design and training techniques for PINNs, including adaptive loss balancing, curriculum training, and a realistic measurement operator.
Regarding blood rheology, the PINN approach accurately estimates viscosity globally and locally under peak systolic conditions. It also provides a clear pattern recognition for diastolic stages. Regarding mass conservation, PINN estimations effectively reproduce the bifurcation of flow through the different branches of the aorta, demonstrate an excellent representation of the non-slip conditions at the walls, and accurately estimate pressure drops with relative errors below the 5 \% in the whole pressure field. We test our pressure drop estimations against the state of the art Virtual Work Energy Relative Pressure (vWERP) estimator, and we observe how our results outperform vWERP in terms of both accuracy and time resolution. Additionally, we find that the best results are achieved by computing the velocity field using the PINN, which is then integrated into the vWERP framework, leading to time super-sampled and high-order approximations, with a clinically admissible accuracy.
Our findings confirm the effectiveness of PINNs in tackling complex inverse problems, as they accurately assimilate observed flow fields, generate precise pressure field estimates, and identify critical regions of high wall shear stress from realistic 4D flow MRI data.
\end{abstract}


\end{frontmatter}

%\tableofcontents

%% main text
\section{Introduction}
%\linenumbers
\label{sec:intro}

%% 4D-flow
4D-flow Magnetic Resonance Imaging (MRI) is an advanced technique that enables the assessment of blood velocities throughout the human body. Unlike 2D-flow MRI or Doppler ultrasonography, 4D-flow provides spatially resolved velocity maps, offering detailed insight into blood flow dynamics~\cite{Azarine2019,Soulat2020,Terada2022}. However, despite its improved spatial resolution, several limitations remain that prevent its use for cardiovascular evaluation. This is where non-invasive, physics-informed estimation techniques become valuable~\cite{Fumagalli2024, Nolte2022}, enabling the inference of high-fidelity velocity, viscosity, and pressure fields from limited or noisy measurements. These reconstructions, in turn, facilitate the computation of clinically relevant biomarkers such as wall shear stress (WSS), oscillatory shear index (OSI), and pressure drops~\cite{Bertoglio2018}.

%% limitations of 4D-flow
Obtaining high-quality 4D-flow MR images usually involves elevated scanning times ($>$20 minutes, depending on the anatomical region and heart rate), even at reduced spatial and temporal resolutions, which often makes the exam not viable from a clinical perspective \cite{Rijnberg2021,Ramaekers2023}. A straightforward approach to reducing scan times is to reduce image resolution.  However, this procedure can lead to significant partial volume artifacts, which compromise the interpretability of the data \cite{Montalba2018}. To accelerate scans, techniques such as parallel imaging, which utilizes the sensitivity of an array of coils, and k-space under-sampling, which takes advantage of the sparsity of the signal in the frequency and time domains, have also been suggested. Unfortunately, both techniques can degrade the Signal-to-Noise Ratio (SNR), ultimately affecting image quality \cite{Demirkiran2022}.

Velocity measurements obtained from 4D-flow MR images are subject to several limitations and potential issues. These include velocity aliasing \cite{Berhane2022}, flow displacement artifacts \cite{Schmidt2021}, and eddy currents \cite{You2022}, all of which can significantly deteriorate the diagnostic quality of the images. Those issues are particularly true when estimating derived quantities such as hemodynamic parameters. Noise and artifacts during the imaging process can lead to considerable discrepancies from the actual physical behavior of blood flow. Additionally, using an inappropriate physics model for blood flow—such as selecting the wrong rheological model—can introduce errors in the quantification of derived quantities \cite{Martinez2023, Yin2025, Mella2026}. Therefore, it is essential to correct 4D-flow velocity measurements for these biases while also incorporating the relevant physical principles for improved estimations.

%4D-flow MR images often face challenges related to limited temporal and spatial resolutions, usually due to hardware constraints. Furthermore, noise, eddy currents, partial volume effects, and flow artifacts can affect the obtained images, making them, in some cases, inadequate to obtain accurate estimates of flow and hemodynamic parameters \cite{Montalba2018,Urbina2016}. These issues underscore the need for improved numerical modeling techniques, like CFD and data-science approaches, to infer missing information and reduce some imaging artifact.

% CFD
Computational Fluid Dynamics (CFD) enables detailed simulations of flow patterns, pressure distribution, and other hemodynamic parameters related to vascular conditions \cite{Carvalho2021}. These simulations have proven to be valuable for pre-surgical planning, studying aneurysm progression, and assessing stenosis \cite{Zhao2014, Li2022, Liu2016, Fumagalli2024}. To effectively integrate 4D flow data with physical laws, optimization strategies are necessary. These strategies ensure that medical measurements closely align with high-fidelity fields that reflect these laws. By incorporating medical measurements, we can address the lack of boundary conditions required to solve the governing equations, allowing us to invert the measurements into a fully resolved 3D state \cite{Galarce2021a, Lombardi2023}. Moreover, CFD simulations often require high computational costs, particularly when traditional numerical methods such as finite volume or finite element methods are applied in three-dimensional dynamic scenarios \cite{Morris2016}.

Recent research highlights the importance of considering shear-thinning blood behavior, even in large vessels \cite{Abugattas2020, Lynch2022, Farias2023}. Furthermore, as discussed in \cite{Martinez2023}, a realistic model of blood rheology has a more significant impact on obtaining biomarkers, such as wall shear stress (WSS), than simply incorporating a turbulence closure model. In \cite{Chauhan2021}, a multi-mode simplified Phan-Thein-Tanner (sPTT) viscoelastic model was calibrated using experimental data, effectively capturing both shear-thinning and viscoelastic effects in simulations of steady and pulsatile flow through stenosed arteries. Similarly, in \cite{Elhanafy2024}, microvascular stenoses were studied using a fully-resolved cellular blood flow model that explicitly simulates red blood cell (RBC) dynamics and platelet margination. Their findings revealed that hematocrit significantly affects cell migration and the formation of a cell-free layer, which, in turn, modulates the apparent viscosity and wall shear stress. This concept is crucial to our work, as we assert that hemodynamic parameters are heavily influenced by blood rheology, which is, in fact, strongly dependent on red blood cell concentration.

In the context of coronary artery biomechanics, the work detailed in \cite{Carpenter2020} emphasizes the critical role of non-Newtonian blood properties in accurately capturing flow features and predicting the initiation of diseases, particularly atherosclerosis. In addition, the study presented in \cite{Qiao2023} investigates the energy dissipation mechanisms in healthy aortas using a fluid-structure interaction (FSI) framework. The findings conclude that viscous friction, significantly influenced by blood rheology, is the primary contributor to energy loss. Together, these studies underscore the importance of incorporating non-Newtonian models and considering hematocrit-related microstructural effects to improve the accuracy of numerical simulations and deepen our understanding of cardiovascular pathophysiology.

%On a different avenue, the definition and modeling of outflow boundary conditions has been a focal point for the hemodynamic community. Researchers have worked on simplified resistance-compliance-resistance strategies to realistically simulate the downstream vascular tree, with Windkessel-based models recognized as the gold standard in this area \cite{Westerhof2009,Piersanti2022,Qiao2019}. In this context, rheological models also have been explored with 

%Recent advances in computational hemodynamics have increasingly emphasized the importance of incorporating the complex rheological behavior of blood and patient-specific physiological features. 

%ROMs, particularly Windkessel models, are widely used to define boundary conditions in vascular simulations. The Windkessel model captures the global properties of arterial systems by representing them as lumped parameter circuits \cite{Westerhof2009}. IN CFD simulations, the three-element RCR (resistance-compliance-resistance) Windkessel model is often used to mimic the behavior of downstream vasculature. By providing a simplified representation of vascular resistance, Windkessel models are essential for generating realistic boundary conditions for blood flow simulations \cite{Piersanti2022,Qiao2019}.

%% inverse problems
Various methodologies have been explored to enhance state estimation accuracy and solve inverse problems. Traditional techniques such as Kalman filters, adjoint methods, and optimization-based approaches have been utilized to infer hidden states or parameters from available data \cite{Huang2011, Fevola2021, Nolte2022}. These methods typically aim to minimize discrepancies between observed data and physical models to estimate unknown quantities, such as pressure fields or material properties \cite{Nolte2022, Lombardi2023}. % However, they may encounter difficulties when dealing with high-dimensional systems and nonlinearities, especially in complex geometries like the cardiovascular system.

Physicians pay particular attention to pressure drops in large arteries, such as the aorta, as these serve as valuable biomarkers for the diagnosis, treatment, and monitoring of conditions like aortic coarctation, valve stenosis, and congenital heart diseases \cite{Warnes2008}. However, the clinical gold standard for assessing such pressure drops remains invasive catheterization, which involves significant risks and high costs for the patient \cite{Zhou2024}. This challenge has sparked increasing interest from the biomedical community, especially given the availability of non-invasive imaging modalities—such as Doppler Ultrasound and 4D-flow MRI—that could potentially be leveraged to estimate pressure fields more safely and cost-effectively.

A classical method for estimating pressure in large vessels is the Pressure Poisson Estimator (PPE) \cite{Ebbers2001}. This technique derives the pressure field by directly substituting velocity measurements into the Navier–Stokes equations, followed by solving a symmetric variational problem for pressure. In \cite{Donati2015}, a related approach is introduced, where the governing equations are tested against the velocity data, leading to a work-energy-based estimate of the pressure drop. An alternative method described in \cite{Svihlova2016} avoids reusing measurements as test functions; instead, it introduces an auxiliary Stokes problem to improve robustness against measurement noise. In \cite{Galarce2021}, a reduced-order model (ROM) approach is employed to reconstruct the whole pressure field, yielding smooth and accurate approximations. This method outperforms traditional techniques, such as the integral momentum-based estimator outlined in \cite{Bertoglio2018}, although it comes with a higher computational cost.
%% SciML

Recently, various machine learning techniques for state estimation have been proposed. However, these methods often face several limitations, including a strong dependence on large datasets, poor generalization, and a lack of interpretability \cite{Alzubaidi2023}. To overcome these challenges, a new subclass of machine learning techniques has emerged, commonly known as Physics-Informed Machine Learning (PIML) \cite{Cuomo2022} or Scientific Machine Learning (SciML) \cite{Raissi2019, Wang2021, Wang2023, Fabra2024}. The core principle of PIML methods is the incorporation of observational, inductive, and learning biases into machine learning models, guiding them toward physically consistent solutions. This approach requires less data while improving both generalization and interpretability \cite{Karniadakis2021, Toscano2025}.

%% PINNs
Among the various techniques in SciML, Physics-Informed Neural Networks (PINNs) have emerged as a promising approach for solving inverse problems, especially in scenarios where the data is limited but the underlying physics is well understood. Unlike traditional neural networks, PINNs incorporate governing physical laws, boundary and initial conditions, and other domain-specific knowledge into their architecture or training processes \cite{Raissi2019}. This integration makes PINNs particularly suitable for medical applications, and a growing body of research has applied PINNs to vascular systems and hemodynamics (see \cite{Toscano2025} for a recent review article). Some studies utilize reduced-order 1D blood flow models, which approximate the vascular tree as a network, allowing significant simplifications of the governing equations. Although this approach sacrifices some accuracy, it enables the efficient modeling of multiple connected vessels \cite{Kissas2020, Sarabian2022, Changdar2024}. More realistic modeling typically occurs in 2D or 3D using simplified geometries and assumes Newtonian rheology for blood \cite{Fathi2020, Arzani2021, DuToit2023}. Other works combine field inference with shape estimation \cite{Sengupta2021}, incorporate fluid-structure interaction \cite{Guan2023, Zhang_Meshless_2024}, focus on turbulent-like flow \cite{Aghaee2025}, or estimate Windkessel parameters \cite{Garay2024}. With respect to the use of PINNs for forward modeling non-Newtonian blood behavior, a study in \cite{Zhang2023} explores a hybrid framework that combines PINNs with PointNet, utilizing ground-truth simulations based on the Carreau-Yasuda constitutive model. In contrast, \cite{Mahmoudabadbozchelou2022} employs PINNs for forward flow simulation in a simplified setting with a wide range of constitutive models.

%% contributions
This work makes a two-fold contribution. First, from a methodological perspective, we utilize advanced techniques such as curriculum training, adaptive collocation points, and adaptive loss balancing \cite{Wang2023,Lu2021,Maddu2022}. These strategies enhance convergence and accuracy, particularly in complex geometries and under challenging data conditions. Second, from an application standpoint, we build upon recent approaches \cite{Garay2024} by refining the observation operator to more effectively replicate the local averaging process inherent in real 4D flow MRI acquisitions \cite{Bissell2023}, thereby narrowing the gap with clinical practice. Additionally, the physical loss in our optimization formulation accounts for the shear-thinning behavior of blood flow \cite{Martinez2023, Mehri2018}, and is explicitly tailored to different hematocrit levels. In summary, we assess the viability of PINNs for estimating viscosity, velocity, and pressure fields from realistic 4D-flow data. Our results suggest that PINNs are helpful for estimating hidden physical fields in the cardiovascular system, providing a potential pathway for improved diagnosis and monitoring of cardiovascular diseases, with particularly accurate pressure drop estimations, relative errors as low as 1\%.

This article is organized as follows. Section \ref{sec:methods} introduces the forward model, describing the boundary value problem for realistic aortic flows under five rheological conditions, and the numerical setup used to solve it. Section \ref{sec:PINN} outlines the PINN-based state estimation methodology, covering model definition, loss function design, and methodological extensions, as well as the generation of synthetic 4D flow data. Section \ref{sec:p_estimation} describes a state-of-the-art method for pressure drop estimation, and a summary of strategies are presented for the estimation of this biomarker. Section \ref{sec:results} presents the reconstructed physical fields, compares them against the reference fields and evaluates key hemodynamic metrics. Finally, Section \ref{sec:conclusion} discusses the quality of the results and potential directions for future research.

\section{Forward Modelling}
\label{sec:methods}

Several strategies exist for modeling blood flow in the cardiovascular system, varying primarily in dimensionality and computational cost. When simulating large arterial trees, or the entire circulatory network, reduced-order models in 0D or 1D are commonly employed. These models rely on spatial averaging of the full 3D Navier–Stokes equations, enabling efficient resolution of global hemodynamic characteristics such as pressure waves and flow distributions with modest computational demands \cite{Montecinos2014}.

In contrast, when localized high-fidelity representations are required, e.g., in stenotic regions or areas of complex flow, the full incompressible Navier–Stokes equations must be solved in three dimensions, capturing detailed spatial variations in velocity and pressure fields \cite{Arbia2016, Deparis2014}. For more realistic modeling, additional complexities, such as fluid-structure interaction (FSI) and electromechanical coupling, can be introduced \cite{Bertoglio2013, Zingaro2024,Ruz2024}. However, in this study, the vessel walls are assumed rigid, an approximation widely adopted for mid- to large arteries such as the aorta and carotids \cite{Galarce2021}.

The rheological behavior of blood is modeled using a non-Newtonian power-law formulation, with parameters explicitly dependent on the hematocrit level. This choice is motivated by experimental \cite{Mehri2018} and numerical evidence \cite{Martinez2023}, which indicate that shear-thinning effects become increasingly relevant under pathological or high-hematocrit conditions. Although some previous studies accept Newtonian assumptions for simplicity \cite{Bertoglio2018, Arbia2016, Galarce2021}, this can result in underestimating key flow characteristics, particularly during diastolic phases or in diseased states \cite{Mella2026}. To address this, we consider a range of hematocrit-dependent viscosity models derived from empirical data measured at physiological temperatures \cite{Walburn1976}, covering both physiological and pathological regimes.

Outflow boundary conditions are implemented using a 3-element Windkessel model at each outlet. These zero-dimensional models simulate the impedance of the downstream vasculature by capturing its resistance and compliance. This approach ensures physiologically consistent pressure wave reflections without explicitly modeling the entire peripheral network \cite{Flores-Geronimo2024, Arbia2016, Formaggia2009}. By combining these elements, we achieve a high-resolution 3D model within a designated area of interest, creating a computationally efficient and physiologically accurate framework.

\begin{figure}[htp]
    \centering
    \includegraphics[width=0.90\textwidth]{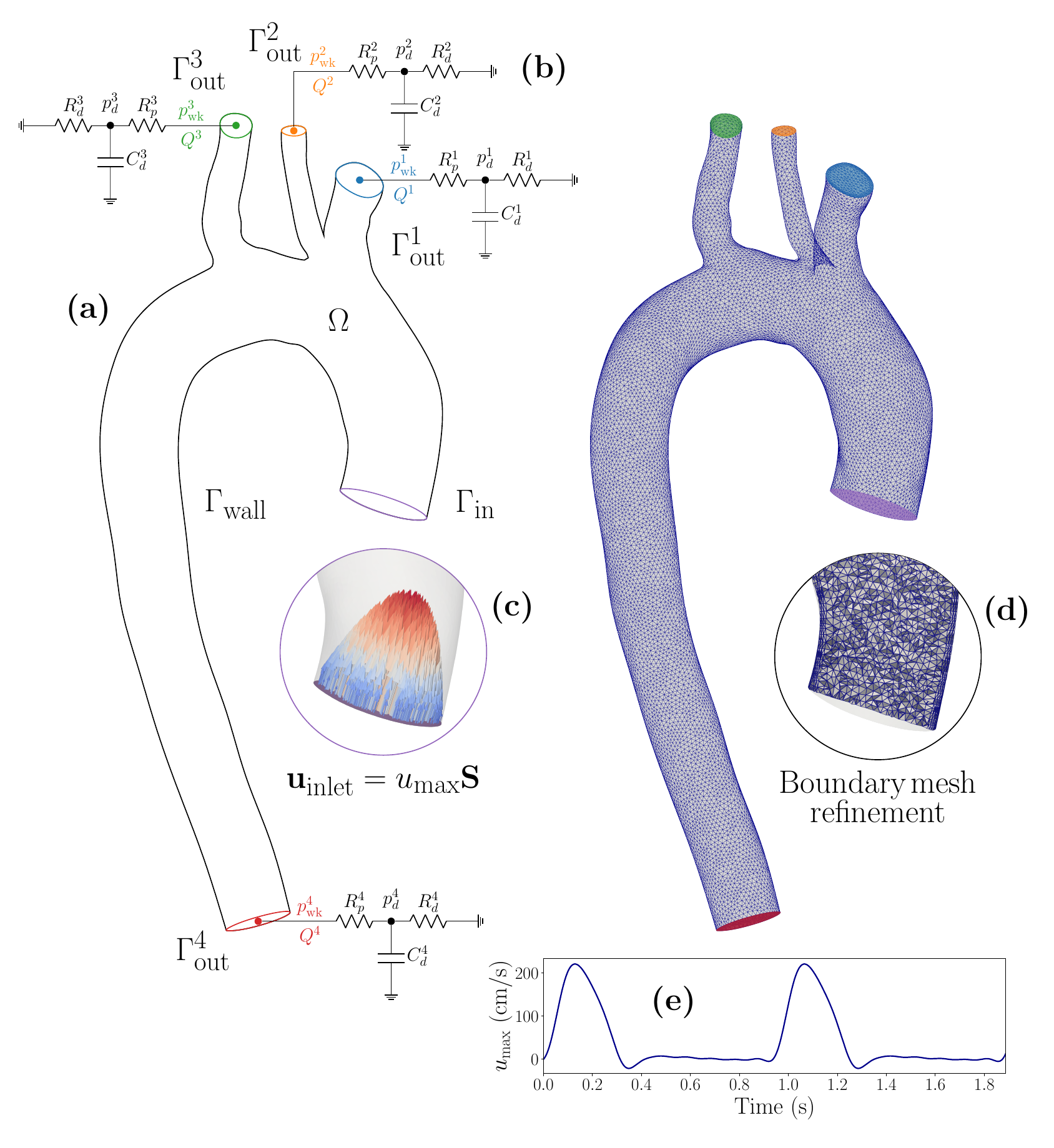}
    \caption{Working set-up: (a) Geometry from the Vascular Model Repository \cite{Updegrove2017}, (b) 0D coupling with Windkessel models for the downstream flow, (c) Inlet paraboloid Dirichlet condition, (d) Tetrahedral mesh refinement near the vessel walls and (e) time-dependent shape of inlet flow.}
    \label{fig:working_domain}
\end{figure}

\subsection{Governing equations}\label{sec:hemoModel}

Let $\Omega\subset\mathbb{R}^3$ be a three-dimensional domain corresponding to the lumen of the thoracic aorta, as shown in Figure \ref{fig:working_domain}. The geometry was obtained from the Vascular Model Repository (Simvascular \cite{Updegrove2017}). The domain boundary $\partial \Omega$ is decomposed into six disjoint sets. That is, the inlet $\Gamma_{\mathrm{in}}$ denoting the connection of the aortic root to the heart; the walls of the vessel $\Gamma_{\mathrm{wall}}$; and the outlets $\Gamma_{\mathrm{out}}^k$, respectively, at the brachiocephalic artery, the left common carotid artery, the left subclavian artery, and the end of the thoracic aorta, such that
\begin{equation}
    \partial\Omega = \Gamma_{\mathrm{in}} \cup \Gamma_{\mathrm{wall}} \cup \bigg( \bigcup_{k=1}^4 \Gamma_{\mathrm{out}}^k \bigg).
\end{equation}

Blood flow is considered an incompressible fluid. The Navier-Stokes equations for our setup adequately describe the phenomena of blood flow and are given by:
\begin{align}
    \rho\left(\frac{\partial\mathbf{u}}{\partial t} + \mathbf{u}\cdot\nabla\mathbf{u}\right) + \nabla p - \nabla \cdot \tau &= 0, &&\text{in} \; \Omega \times [0,T], \label{eq:moment}  \\
    \nabla \cdot \mathbf{u} &= 0, &&\text{in} \; \Omega \times [0,T], \label{eq:mass}\\
    \mathbf{u} &= \mathbf{u}_{\mathrm{inlet}}, &&\text{on} \; \Gamma_{\mathrm{in}} \times [0,T], \label{eq:inletflow} \\
    \mathbf{u} &= \mathbf{0}, &&\text{on} \; \Gamma_{\mathrm{wall}} \times [0,T], \label{eq:noslip} \\
    (\tau - p I) \cdot \mathbf{n} &= -p_\mathbf{wk}^k\mathbf{n},&&\text{on} \; \Gamma_\mathrm{out}^k \times [0,T], \label{eq:compat} 
\end{align}
where $\mathbf{u} : \Omega \rightarrow \bR^3$ represents the velocity vector field and $p : \Omega \rightarrow \bR$ the pressure field. The fluid density is denoted by $\rho$, and 
$\tau$ represents the shear stress tensor which includes non-Newtonian rheology. The remaining terms relate to boundary conditions, which are discussed next. %The inlet boundary condition, $\mathbf{u}$, is a pulsatile parabolic inflow adapted to the cardiac frequency and the mass flow rate. 

\subsection{Boundary conditions}

The inlet flow condition is a function $\mathbf{u}_\mathbf{inlet}(\mathbf{x},t) = u_\mathbf{max}(t) \mathbf{S}(\mathbf{x})$, where $u_\mathrm{max}(t)$ is the maximum flow velocity over time obtained from MR measurements, and $\mathbf{S}(\mathbf{x})$ is a normalized paraboloid function \cite{LaDisa2011}, as illustrated in Figure \ref{fig:working_domain}. The parameters for the Windkessel model at each outlet are summarized in Table \ref{tab:wkparams}. These values are also an input for this work and were obtained from the Vascular Model Repository.

\begin{table}[h]
  \centering
  \begin{tabularx}{0.9\textwidth}{>{\centering\arraybackslash}X>{\centering\arraybackslash}X>{\centering\arraybackslash}X>{\centering\arraybackslash}X>{\centering\arraybackslash}X}
    \hline
    Parameter & $\Gamma_\mathrm{out}^1$ & $\Gamma_\mathrm{out}^2$ & $\Gamma_\mathrm{out}^3$ & $\Gamma_\mathrm{out}^4$ \\[1ex]
    \hline
        Proximal resistance $R_p^k$ &  274 & 1300 & 791 & 141 \\
        Distal resistance $R_d^k$ &  5675 & 19663 & 10048 & 2066 \\
        Capacitance $C_d^k$ &  5.08 & 1.4416 & 2.788 & 13.6904 \\
        Initial distal pressure $p_d^k$ & 107325 & 107325 & 107325 & 107325 \\
    \hline
  \end{tabularx}
  \caption{Parameters for the outlet Windkessel models.}\label{tab:wkparams}
\end{table}

A no-slip boundary condition is applied to the vessel walls. On the outlets, a three-element Windkessel model is used to represent the downstream flow:
\begin{align}
    &Q^k = \int_{\Gamma_\mathrm{out}^k} \mathbf{u} \cdot \mathbf{n} \;\mathrm{dA} \label{eq;WKflowrate},\\
    &p_\mathrm{windk}^k = R_\mathrm{p}^k Q^k+ p_\mathrm{d}^k \label{eq:WKpressure},\\
    &C^k \frac{dp_d^k}{dt} + \frac{p_\mathrm{d}^k}{R_\mathrm{d}^k} = Q^k. \label{eq:WKcontinuity}
\end{align}
A Windkessel model is a lumped parameter model that captures the effects of downstream vasculature on the outflow, via coupling with a set of ordinary differential equations of a distal pressure $p_d^k$ at each outlet. Consequently, $R_\mathrm{p}^k$, $R_\mathrm{d}^k$ represent resistance to the proximal and distal flow of the outlet, while capacitance $C^k$ represents the compliance of the distal vessels. Equations \eqref{eq:WKpressure} and \eqref{eq:WKcontinuity} define the dynamic response of $p_d^k$ to the uniform normal stress at the outlet $p_\mathrm{wk}^k$ and the outflow rate $Q^k$. In turn, $p_\mathrm{wk}^k$ affects the velocity and pressure fields via the compatibility relation in equation \eqref{eq:compat} \cite{Formaggia2009}. This modeling approach follows previous work by the authors \cite{Galarce2021a, Galarce2021}.

\subsection{Blood rheology}

The shear thinning behavior is modeled using a power-law constitutive relation. That is, the shear stress tensor $\tau$ is related to the strain rate tensor $\dot\varepsilon = \frac{1}{2}(\nabla\mathbf{u} + \nabla\mathbf{u}^\intercal)$ by an apparent viscosity $\mu_\mathrm{PL}$ as follows:
\begin{equation}\label{eq:gennewtonian}
    \tau = 2 \mu_\mathrm{PL} \dot\varepsilon,
\end{equation}
\begin{equation}\label{eq:MUPL}
    \mu_\mathrm{PL} = m \dot\gamma^{n-1},
\end{equation}
where $m$  is the consistency index in ($\mathrm{Pa}\,\mathrm{s}^n$), $n$ the power-law index, and $\dot\gamma = \sqrt{2\dot\varepsilon : \dot\varepsilon}$ the strain rate. The power-law index is constrained to $0<n\leq 1$ to model shear thinning behavior with $n=1$ recovering the Newtonian case \cite{Abugattas2020}.

Since the article focuses on the estimation of the state of blood flow that incorporates realistic rheology, five CFD experiments were performed, which encompassed five different levels of hematocrit-ranging from anemic at 20\% to polycythemic at 70\%. All cases are modeled with the power law relation defined in Section \ref{sec:hemoModel}. Table \ref{tab:rheoparams} shows all the relevant parameters for each simulation.

\begin{table}[ht]
  \centering
  \begin{tabularx}{0.6\textwidth}{>{\centering\arraybackslash}X>{\centering\arraybackslash}X>{\centering\arraybackslash}X}
    \hline
    Hematocrit ($\%$) & Power-law \, $m \, (10^{-2} \mathrm{Pa}\,\mathrm{s}^n)$, $n$ \\[1ex]
    \hline
        20.0 & 0.6850, 0.7113 \\
        32.5 & 1.7271, 0.6339 \\
        45.0 & 2.4208, 0.7146 \\
        57.5 & 4.1933, 0.6349 \\
        70.0 & 5.3985, 0.6313 \\
    \hline
  \end{tabularx}
  \caption{Rheological parameters for the CFD experiments.}\label{tab:rheoparams}
\end{table}

\subsection{Numerical set-up}

To generate a set of ground truth solutions, we conduct computational fluid dynamics (CFD) simulations using the finite element method (FEM) with equal-order linear interpolations for both velocity and pressure. To circumvent the inf-sup saddle point condition, we apply standard SUPG–PSPG stabilization techniques \cite{Volker2016}. For time discretization, we use a second-order backward differentiation formula, while Windkessel boundary conditions are solved using a fourth-order Runge-Kutta scheme. 

The computational domain is discretized in space using a tetrahedral mesh of 127131 vertices and 719419 elements, using an adaptive wall refinement to better capture the effects of the boundary layer. For temporal discretization, we set a time step of $\Delta t = 9.37 \cdot 10^{-4} \text{ s}$ to cover two complete cardiac cycles, leading to 1000 steps per cycle. We consider simulations only once the periodic cardiovascular regime is reached. The non-linear convective term in the momentum equations is managed using a Picard scheme, with a relative residual $H^1(\Omega) \times L^2(\Omega)$ tolerance of $10^{-4}$.

The forward model will be used twofold. On the one hand, it will be used to design the optimization problem cost function in the next section. On the other hand, the forward model will be used to generate synthetic data to use as ground truth solutions. Our work will produce those simulations using the finite element method with the software MAD \cite{galarceThesis, galarceMAD}, which works upon the CPU parallel linear algebra PETSc solver \cite{petsc}.

\section{State estimation with Physics Informed Neural Networks}\label{sec:PINN}

This section outlines the process of conducting state estimation using a Physics-Informed Neural Network (PINN) in conjunction with 4D-flow MRI data. We introduce a neural field representation for both velocity and pressure fields, and discuss the design of the loss function, the training procedure, extensions of the PINN, the generation of synthetic 4D flow data, and details regarding the implementation of neural network training.

We employ a fully connected feed-forward neural network as a universal approximator of the problem outlined in Section \ref{sec:hemoModel}, including the no-slip boundary conditions and the 4D-flow data in the optimization process. Let $\mathrm{NN}: \mathbb{R}^4 \rightarrow \mathbb{R}^4 $ be a neural field that maps a point in space and time $(\mathbf{x},t)$ to an estimation of physical fields $(\hat u, \hat v, \hat w, \hat p)$, such that:
\begin{equation}\label{eq:NN}
    \mathrm{NN}(\mathbf{x},t;\theta) = \left( f_L \circ f_{L-1} \circ \ldots \circ f_1 \right) (\mathbf{x}, t) = \begin{bmatrix} \hat u & \hat v & \hat w & \hat p \end{bmatrix},
\end{equation}
where each layer is defined with the affine transformation $f_k = W^{(k)}\sigma(f_{k-1}) + b^{(k)}$ and the set of all parameters is denoted $\theta = \{W^{(\ell)}, b^{(\ell)}\}_{\ell=0}^L$. We can find the approximate solution by solving a multi-objective optimization problem over $\theta$, minimizing an appropriately designed loss function that accounts for the residual minimization criteria of the underlying physical laws and the discrepancy between the data and the estimation when projected onto low-resolution voxels. We can solve the optimization problem using stochastic gradient descent. Furthermore, with a suitable choice of activation function $\sigma$, this representation is differentiable with respect to both its parameters and input, which allows the computation of the first and second space and time derivatives required by the physical laws residuals. 

\subsection{Loss Function Design}\label{sec:loss}

We approach the multi-objective optimization problem by designing a total loss that adds together contributions from each objective as:
\begin{equation}
    \mathcal{L}(\theta) = \mathcal{L}_{obs}(\mathcal{D}_{obs}, \theta) + \mathcal{L}_{p}(\mathcal{D}_{p}, \theta) + \mathcal{L}_{NS}(\mathcal{D}_\mathrm{int}, \theta) + \mathcal{L}_{BC}(\mathcal{D}_\mathrm{bd}, \theta).
\end{equation}
Each term is a mean square loss evaluated over a different set of points $\mathcal{D}_{(i)}$ with weightings $\lambda_{(i)}$ to modulate its contribution to the total loss. We next break down each loss component:
\begin{enumerate}[leftmargin=0.5cm]
    %%%% 4D FLOW
    \item The observation loss $\mathcal{L}_{obs}$ quantifies the error concerning the 4D-flow measurements obtained from MR images. The flow field is biased for low-resolution images due to imaging artifacts (especially partial volume effects). We introduce an observer matrix $W$ \cite{Brunton2022} that mimics the observation process, taking local averages within the voxels, from simulated field values, thus leading to an unbiased sub-sampling of the high-fidelity solutions. This approach is a mid-point between realism and other approaches, such as direct nodal evaluation \cite{Garay2024}, or sub-sampling by interpolation on coarser meshes \cite{Bertoglio2018}. We compute the local averages using a quasi-Monte Carlo (QMC) integration scheme. The quadrature points are drawn from the Halton sequence \cite{Halton1960} and randomized with Owen scrambling \cite{Owen2017}. Although only a modest number of quadrature points are required per voxel, the total number needed for batch optimization is prohibitively large. Thus, to implement mini-batching under this setting, the quadrature points are drawn, scaled, and translated one random mini-batch of voxels at a time. This also provides the regularization benefits of mini-batch optimization. In summary, we compute the error by comparing the values of a mini-batch of voxels from the images with the corresponding values generated from the model:
    \begin{equation}
        \mathcal{L}_\mathrm{obs} =
        \lambda_\mathrm{obs} \frac{1}{N_\mathrm{obs}} \sum_{i=0}^{N_{obs}} \bigg(\mathbf{u}^{(i)}_\mathrm{obs} - W \mathbf{\hat{u}}(\mathbf{x}^{(i)};\theta) \bigg)^2.
    \end{equation}
    
    %%%% GLOBAL MEAN PRESSURE
    \item The pressure data loss, denoted as $\mathcal{L}_{p}$, represents the error associated with the measurement of the average physiological pressure over the entire space-time domain. Since 4D-flow techniques only measure velocities, the pressure field must be fully inferred. This fact presents a challenge, as the pressure can oscillate during training, potentially leading to instabilities. To address this issue, we introduce minimal pressure data and utilize an averaging operator $W_p$ that computes local averages, similar to how the velocity observer $W$ operates. This approach simplifies this component of the loss to:
    \begin{equation}\label{eq:Lp}
        \mathcal{L}_\mathrm{p} =
        \lambda_\mathrm{p} \frac{1}{N_\mathrm{p}} \sum_{i=0}^{N_{p}} \bigg(p_\mathrm{mean} - W_p \hat p(\mathbf{x}^{(i)};\theta) \bigg)^2.
    \end{equation}

    %%%% NAVIER-STOKES + CONTINUITY
    %\item The physical loss $\mathcal{L}_\mathrm{NS}$, adds the residuals for the conservation of momentum and mass equations. The required derivatives are computed using automatic differentiation (AD).
    
    %Since evaluating the residuals at a small, fixed set of points can lead to model overfitting, we compute them over approximately 20 million pre-computed, scrambled Halton points, leading to a different cloud per iteration. Thus, the component-wise physical cost function reads:

\item The physical loss term, denoted by $\mathcal{L}_\mathrm{NS}$, aggregates the residuals of the momentum and mass conservation equations. These residuals are evaluated using derivatives computed via automatic differentiation (AD).

To prevent overfitting that may arise from evaluating residuals at a small, fixed set of points, we instead compute them over approximately 20 million precomputed Halton points, which are scrambled at each iteration. This ensures that the residuals are evaluated over a different point cloud in every training step. Accordingly, the component-wise physical loss function is given by:
    \begin{equation}\label{eq:LNS}
        \begin{aligned}
            \mathcal{L}_\mathrm{NS} &=
            \lambda_\mathrm{NSx} \frac{1}{N_\mathrm{int}} \sum_{i=0}^{N_\mathrm{int}} \bigg( \mathcal{NS}[\hat u, \hat p] \bigg)^2
            +
            \lambda_\mathrm{NSy} \frac{1}{N_\mathrm{int}} \sum_{i=0}^{N_\mathrm{int}} \bigg( \mathcal{NS}[\hat v, \hat p] \bigg)^2 \\
            &+ 
            \lambda_\mathrm{NSz} \frac{1}{N_\mathrm{int}} \sum_{i=0}^{N_\mathrm{int}} \bigg( \mathcal{NS}[\hat w, \hat p] \bigg)^2
            +
            \lambda_\mathrm{Cont} \frac{1}{N_\mathrm{int}} \sum_{i=0}^{N_\mathrm{int}} \bigg(\nabla \cdot \mathbf{\hat{u}}(\mathbf{x}_\mathrm{int}^{(i)};\theta)\bigg)^2,
        \end{aligned}
    \end{equation}
    where $\mathcal{NS}$ is the operator of the momentum conservation equations, which is defined as follows:
    \begin{equation}\label{eq:NS}
        \mathcal{NS}[\hat u_i, \hat p] = \rho \bigg(\frac{\partial\hat u_i}{\partial t} + (\mathbf{\hat u} \cdot \nabla) \hat u_i \bigg) + \frac{\partial\hat p }{\partial x_i} - (\nabla \cdot \tau(\hat u_i))_i.
    \end{equation}

    This loss is evaluated point-wise in the interior of the domain, computing the required derivatives using automatic differentiation. Since a small and fixed set of collocation points can lead to model overfitting, we sample a set of approximately 20 million scrambled Halton points. Dense sampling promotes the learning of physically consistent solutions across the domain.

    %%%% NO-SLIP CONDITION
    \item We incorporate an additional term, denoted as $\mathcal{L}_\mathrm{BC}$, that addresses the no-slip boundary conditions on the vessel walls. This term is calculated within two specific sub-domains: i) the vessel wall itself, and ii) a boundary volume outside the vessel wall, which has a thickness of 1.5 times the largest voxel dimension. The inclusion of this outer boundary encourages the network to assign zero velocity to points outside the domain. The chosen thickness ensures that when we compute local averages of voxels adjacent to the vessel wall, any quadrature point outside the working domain is set to have zero velocity. The volume boundary points consist of 14 million precomputed, randomized Halton points. During the training phase, mini-batches are drawn from this set, maintaining a fixed ratio of mesh points to boundary volume points. The corresponding loss function is defined as follows:
    \begin{equation}
        \mathcal{L}_\mathrm{BC} =
        \lambda_\mathrm{BC} \frac{1}{N_\mathrm{BC}} \sum_{i=0}^{N_\mathrm{BC}} \bigg(\mathbf{\hat u}(\mathbf{x}_\mathrm{BC}^{(i)};\theta)\bigg)^2.
    \end{equation}
    %, as a compromise between the cost of enforcing zero velocity for a large volume outside the domain and ensuring that when computing the local velocity average for a voxel at the edge of the domain, quadrature points that fall outside have zero velocity. 
    
\end{enumerate}

\begin{figure}[!ht]
    \centering
    \includegraphics[width=1\linewidth]{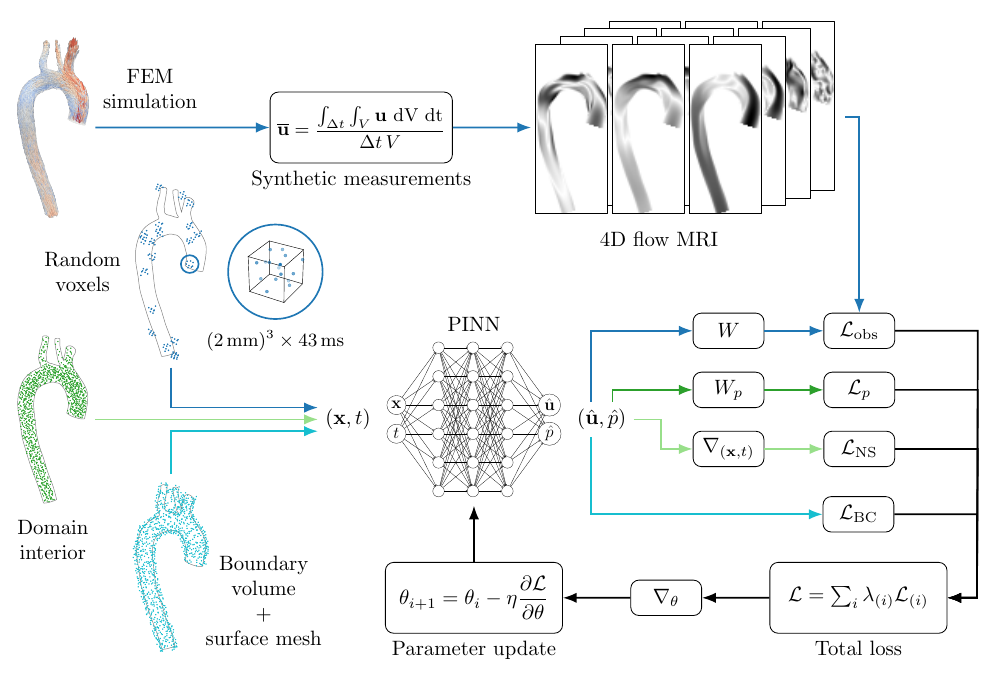}
    \caption{Computation pipeline for a single training step. The observation loss is fed with 4D-flow like local averages in time and space of the velocity field computed with quasi-random quadrature points. Once the total loss is assembled, the learning optimization process takes place to iteratively approach the optimal parameter set $\theta$.}
    \label{fig:diagram}
\end{figure}

%%% PINNS EXTENSIONS %%%
\subsection{PINNs Extensions}\label{sec:pinnextensions}

We address a three-dimensional, time-dependent problem governed by a non-linear system of partial differential equations, with supervision provided only through indirect data. This setting presents significant challenges for standard PINN approaches. To overcome these difficulties, we incorporate several additional techniques to enhance performance, namely:

\begin{enumerate}
    %%%% DATA SCALING/ADIM
    \item A key machine learning technique for effective training involves scaling the data to normalize or standardize it. Accordingly, we introduce the following rescaled variables: 
    When modeling physical variables, we also apply non-dimensionalization for the same purpose. 
    \begin{align*}
        x_i &= \frac{x_i - \tilde{x}_i}{L}, &
        t   &= f\cdot(t - t_\mathrm{min}), &
        u_i &= \frac{u_i}{U}, &
        p   &= \frac{p}{p_0},
    \end{align*}
    where $L$ is the mean length of the mesh bounding box; $\tilde{x}_i$ are shifts on each space coordinate that center the domain inside the voxel grid, which itself begins at the origin; $f$ is the heart rate expressed in Hertz; $t_\mathrm{min}$ is the time at which the cardiac cycle begins; $U$ and $p$ are characteristics values for velocity and pressure near the peak systolic regime, of $200\,\mathrm{cm/s}$ and $1.86\cdot 10^{5}\, \mathrm{Ba}$, respectively.
    % no se si poner una tabla para estos 6 numeritos
    % L = 12.06 cm
    % xs = -8.7628 cm
    % ys = -2.9242 cm
    % zs = -20.1243 cm
    % f = 1 / 0.93 s = 1.068 1/s
    % tmin = 0 para las imagenes, = 0.93 para el ideal test set.

    %%%% QUASI-MONTE CARLO + AGGRESIVE RESAMPLING
    \item Throughout this study, several integrals need to be evaluated efficiently, in an easily differentiable way. As has already been noted, we approach this problem using quasi-Monte Carlo integration. To limit GPU memory use, we use a minimal number of quadrature points which are re-drawn at each training step. Trivial shape domains, like a voxel, can be sampled during training, while interior and boundary quasi-random points are pre-computed via rejection sampling using the surface of the vessel. We prefer quasi-random points over pseudo-random points because they exhibit similar or even better convergence properties for Monte Carlo integration and are faster to sample \cite{Halton1960}.

    %%%% NETWORK PRECONDITIONING
    \item To facilitate the learning of a complex field, we employ a two stage curriculum training scheme \cite{Wang2023, Toscano2025}. First, in the preconditioning stage, the model is trained on the 4D-flow field, interpreted as point measurements on the center of each voxel and without any physics loss terms. Second, the physics training stage is performed with all losses as described in Section \ref{sec:loss}. The first stage acts as an initialization for the neural network. This decision is motivated by the fact that fitting MR images is fast and gives a good initial estimate of the flow.

    %%%% LOSS BALANCING
    \item  We use an adaptive weighting technique proposed in \cite{Maddu2022}. This method balances the loss components by equalizing the variance of the gradients of each component with respect to $\theta$. This makes training more robust against catastrophic forgetting. Otherwise, when the model begins the physics training stage, the gradients of the not-yet minimized losses would overwhelm the vanishing gradient of the already small observations loss, eliminating the progress made in the preconditioning stage. The weightings update rule follows:
    \begin{align}
        \hat \lambda_{(i)}(\tau) &= \Bigg( {\frac{\mathrm{std}\{\nabla_\theta \mathcal{L}_{(i)}(\tau) \}}{\underset{i=1,\dots,K}{\mathrm{max}} \big(\mathrm{std}\{\nabla_\theta \mathcal{L}_{(i)}(\tau) \}\big)}} \Bigg)^{-1},
         \label{eq:maddu1} \\
        \lambda_{(i)}(\tau+1) &= \alpha\lambda_{(i)}(\tau) + (1-\alpha)\hat \lambda_{(i)}(\tau), \label{eq:maddu2}
    \end{align}
    where $\mathcal{L}_{(i)}$ represent each loss component, $\lambda_{(i)}$ each loss weighting, and $\nabla_\theta$ the gradient with respect to the network parameters. Equation \eqref{eq:maddu1} defines the $i$th new weight as the reciprocal of the standard deviation of the components of $\nabla_\theta \mathcal{L}_{(i)}$, normalized by the maximum standard deviation among all the losses. Then, to prevent instability due to the changing weights, especially if the update is not computed in every step, Equation \eqref{eq:maddu2} applies an exponential moving average.

    %%%% ADAPTIVE COPOINT REFINEMENT
    \item To improve physics learning, we use Residual-based Adaptive Refinement with Greed (RAR-G) for collocation points refinement \cite{Lu2021,Wu2023}. RAR-G is used on a parallel set of points, in addition to the large set of precomputed quasi-random points described in Section \ref{sec:loss}. The samples for RAR-G updates are taken from each batch of collocation points used to compute $\mathcal{L}_{\mathrm{NS}}$ at each iteration.

\end{enumerate}

The parameter estimation of this work is carried out on a 16 GB VRAM Nvidia RTX A4000 GPU. We mimic the network architecture of \cite{Garay2024}, where it is verified that the neural field holds sufficient expressivity for aortic flows. The 4D-flow data is organized as a design matrix, and randomly split into a training set with $90\%$ of the observations, a validation set with $5\%$, and a test set with the remaining $5\%$. As described in Section \ref{sec:pinnextensions}, a two-stage scheme is followed. For both stages of training, the \textit{Adam} optimizer \cite{Kingma2014} is used together with a \textit{Reduce Learning Rate on Plateau} scheduler. The scheduler tracks a validation loss, which is computed in the same way as the observation loss but using the validation set. All code is written in Python using the PyTorch library \cite{Paszke2019}. The optimization hyperparameters are shown in Table \ref{tab:trainingparams}. Preconditioning training runs take $13$ to $14$ minutes each, and physics training runs take $\sim7.5$ hours each. Training times are enough to complete the learning stage, so losses reach a plateau condition.

\subsection{Synthetic Training Data}\label{sec:traindata}

4D-flow MRI is a technique that captures blood flow velocities within a three-dimensional volume throughout the cardiac cycle. This imaging process utilizes bipolar magnetic gradient pulses to encode velocity information in all three spatial dimensions. The phase difference observed between two acquisitions is directly correlated with the velocity of the moving spins. Additionally, these acquisitions are synchronized with different phases of the cardiac cycle using electrocardiogram (ECG) gating \cite{Soulat2020}. It is important to note that this work assumes that voxels in MRI represent local spatiotemporal averages of velocity, without taking into account noise or the more complex artifacts typically associated with MR imaging.

Synthetic images were generated from CFD simulations. The local averages are computed by linearly interpolating the FEM solution to a grid and then averaging the grid points inside each voxel. Points outside the convex hull of the mesh are computed using nearest neighbors. The mesh is used to generate an exact segmentation of the flow.

All sets of images have a voxel size of $2\times2\times2\, \mathrm{mm}^3$, dividing the spatial domain into a $24\times47\times111$ grid and the cycle into 22 phases with a duration of $42.6$ ms. According to the 2023 4D-flow consensus \cite{Bissell2023}, these parameters correspond to low-quality observations. Thus, the synthetic velocity measurement at a particular voxel reads:
\begin{equation}
    \overline{\mathbf{u}} = \frac{1}{(t_{i+1} - t_i) \, V} \int_{t_i}^{t_{i+1}} \int_\mathrm{voxel} \mathbf{u} \text{ dV dt},
    \label{eq:voxels}
\end{equation}
where $V$ refers to the voxel volume, while $t_i$, $t_{i+1}$ are instants of time of the start and end of the signal acquisition interval defining a cardiac phase, such that $t_{i+1} - t_i = 42.6$ ms. Measurements of $\mathbf{u}$ are acquired for all voxels in a grid and all cardiac phases.

Our approach offers a novel perspective compared to that of \cite{Garay2024}, which models the imaging process solely as point measurements taken at control points directly derived from the forward simulation. In contrast, \cite{Bertoglio2018} creates synthetic measurements by interpolating the high-fidelity solution onto a lower-resolution mesh. While the work in \cite{Galarce2021} does employ a spatial averaging strategy similar to ours, it is particularly noteworthy that it overlooks the intrinsic time delay introduced by the progressive acquisition of magnetic resonance imaging in the frequency domain. This time delay consideration represents an advancement in our methodology, enhancing the accuracy and realism of the imaging process.

The mean pressure in Equation \eqref{eq:Lp} is computed as an average over the entire space-time domain, $\Omega\times[0,T]$::
\begin{equation}
    p_\mathrm{mean} = \frac{1}{T \, \vert\Omega\vert} \int_{0}^{T} \int_\Omega p  \text{ dV dt}.
\end{equation}

The results of the mean pressure for each hematocrit level are shown in Table \ref{tab:pmeans}.
\begin{table}[h]
  \centering
  \begin{tabularx}{0.4\textwidth}{cc}
    \hline
    Hematocrit ($\%$) & $p_\mathrm{mean}\ (\mathrm{Ba})$ \\[1ex]
    \hline
        20.0 & 136381 \\
        32.5 & 136401 \\
        45.0 & 136489 \\
        57.5 & 136491 \\
        70.0 & 136530 \\
    \hline
  \end{tabularx}
  \caption{Computed global mean pressures for each experiment.}\label{tab:pmeans}
\end{table}

\begin{table}[!hb]
\parbox{.4\linewidth}{
\begin{tabular}{ll}
    \hline
    \multicolumn{2}{c}{Network Architecture} \\
    \hline
        Hidden layers       &  6             \\
        Neurons per layer   &  256           \\
        Activation function &  Swish (SiLU)  \\
        Param. init.        & Kaiming normal \\
        \hline
        & \\
        & \\
    \hline
    \multicolumn{2}{c}{Preconditioning Stage} \\
    \hline
        Epochs                & 200        \\
        Batch size            & 1024       \\
        Optimizer             & Adam       \\
        Initial learning rate & 0.001      \\
        Scheduler             & Reduce LR  \\
                              & on Plateau \\
        Factor                & 0.5        \\
        Patience              & 10         \\
    \hline
\end{tabular}
}
\hspace{1cm}
\parbox{.4\linewidth}{
\begin{tabular}{ll}
    \hline
    \multicolumn{2}{c}{Physics Training Stage} \\
    \hline
        Epochs                & 200               \\
        Batch observations    & 512               \\
        Points per voxel      & 16                \\
        Batch size interior   & 4096              \\
        Batch size refinement & 2048              \\
        Refinement            & 128               \\
        re-sample size        &                   \\
        Batch size boundary   & 4096              \\
        condition             &                   \\
        Optimizer             & Adam              \\
        Initial learning rate & 0.001             \\
        Scheduler             & Reduce LR         \\
                              & on Plateau        \\
        Factor                & 0.5               \\
        Patience              & 5                 \\
        Loss balancing        & Inverse Dirichlet \\
        alpha                 & 0.99              \\
        Update period (steps) & 10                \\
    \hline
\end{tabular}
}
\caption{Network design and training hyperparameters.}\label{tab:trainingparams}
\end{table}

\section{Relative pressure estimation}\label{sec:p_estimation}

A relevant quantity of interest for the biomedical community is the pressure drop between a given inlet-outlet surface computational domain:
\begin{equation} \label{eq:pdelta}
\delta p_{k} = \frac{1}{|\Gamma_\mathrm{in}|} \int_{\Gamma_\mathrm{in}} p \text{ dx} - \frac{1}{|\Gamma_\mathrm{out}^k|} \int_{\Gamma_\mathrm{out}^k} p \text{ dx}, \quad k=1,2,3,4.
\end{equation}

Several correlations have been made over the years to link pressure drop with heart function \cite{Courtois1990, Pasipoularides1987}. We next briefly describe the state-of-the-art in pressure drop estimation and set a benchmark framework to assess the quality of the PINN reconstruction against the best current strategy.

\subsection{State-of-the art methods}

Among the various alternatives available in the literature to compute pressure drops, a classical option is the pressure Poisson estimator \cite{Ebbers2001}. This estimator has been tested in hemodynamic problems alongside newer techniques such as the Stokes estimator \cite{Svihlova2016} and the work-energy relative pressure estimator \cite{Donati2015}. Recent research highlights the advantages of the so-called virtual Work-Energy Relative Pressure (vWERP) estimator \cite{Bertoglio2018, Marlevi2021}, particularly in terms of its accuracy and robustness against noise. In the following section, we will describe the vWERP and conduct numerical experiments using synthetic data to assess the quality of pressure drop reconstruction compared to the PINN approach.

The vWERP method works upon the assumption that the 4D-flow measurements satisfy the weak form of the Navier-Stokes equations. The problem is thus reformulated for a synthetic measured field $\bu_m$ as follows:
\begin{align}
 - \int_{\Omega} \nabla p~ \xi \dx &= \rho \int_{\Omega} \frac{\partial{\bu_m}}{\partial t} \xi ~\dx + \rho \int_{\Omega} \bu_m \cdot \nabla \bu_m \xi ~\dx \\
 &+ \int_{\Omega} 2\mu_{\text{PL}} \nabla^s(\bu_m) : \nabla^s(\xi) \dx - \int_{\partial \Omega} 2 \mu_{\text{PL}} \nabla^s \bu_m \cdot \bn ~\xi ~\ds,
\end{align}
$\forall \xi \in V$, where $V$ is an ad-hoc Sobolev space. We remark that another underlying assumption of the method is related to the linear rheology of the fluid, yet we consider that introducing the shear-thinning behavior in the method is a mere modification that does not justify a new name for it.

The central aspect of the methodology is to effectively manage the pressure term. By carefully selecting the test functions $\xi$, we can directly calculate the pressure drop at each outlet of the domain. To begin this process, we will first apply integration by parts:
\begin{equation}
    -\int_{\Omega} \nabla p ~\xi \dx = \int_\Omega p \nabla \cdot \xi \dx - \int_{\partial \Omega} p ~\xi \cdot \bn \ds,
\end{equation}
and we notice that, to recover mean pressures from $u_m$, one can select fixed $\xi_k^*$ that is solenoidal ($\nabla \cdot \xi_k^* = 0$), by solving an auxiliary Stokes problem for each outlet: Find $\xi_k^* \in V$ such that:
\begin{align}
-\Delta \xi_k^* + \nabla p^* &= 0 &&\text{in } \Omega, \\
\nabla \cdot \xi_k^* &= 0 &&\text{in } \Omega, \\
\xi_k^* \cdot \bn &= -\bn &&\text{on } \Gamma_\text{in}, \\ 
(\nabla^T + \nabla) \xi_k^* - p\bn &= 0 &&\text{on } \Gamma_\text{out}^k, \\ 
\xi_k^* &= 0 &&\text{on } \partial \Omega \setminus  (\Gamma_{\text{out}}^k \cup \Gamma_\text{wall}),
\end{align}
for $k=1,2,3,4$. In addition, if we assume constant pressure in the boundary inlets and outlets, we can derive the following pressure estimator $\delta p_k$ (for the outlets $k=1,2,3,4$):
\begin{align}
\delta p_k &=
- \frac{1}{A_{\text{inlet}}} \Bigg(
\underbrace{\rho \int_{\Omega} \frac{\partial{\bu_m}}{\partial t} \xi_k^* ~\dx}_\text{virtual kinetic power} +
\underbrace{\rho \int_{\Omega} \bu_m \cdot \nabla \bu_m \xi_k^* ~\dx}_\text{virtual convective power} \\
&+
\underbrace{\int_{\Omega} 2 \mu_{\text{PL}} \nabla^s(\bu_m) : \nabla^s(\xi_k^*) \dx - \int_{\partial \Omega} 2 \mu_{\text{PL}} \nabla^s \bu_m \cdot \bn \xi_k^*~\ds}_\text{virtual viscous power}
\Bigg),
\label{eq:vWERPpdrop}
\end{align}
where $A_{\text{inlet}} = - \int_{\Gamma_{\text{inlet}}} \xi^*_k \cdot \bn \ds = \int_{\Gamma_{\text{inlet}}} \bn \cdot \bn \ds$ is the inlet area. 

\subsection{Numerical implementation and summary of strategies}

We implement the vWERP strategy using $\bP_1$ linear piece-wise Lagrange elements in the MAD \cite{galarceThesis, galarceMAD} finite element modular library. Concerning the time scheme, we test the method with a first- and second-order finite-difference formula for the kinetic energy time derivative. The remaining terms are evaluated with a mid-point time rule.

In summary, we will benchmark the pressure drop reconstruction with the following approaches:

\begin{enumerate}
\item PINN: we rely on joint estimation of the neural field from 4D-flow measurements, extracting only the pressure field and calculating the pressure drop using the appropriate formula \eqref{eq:pdelta}. 
\item vWERP: we test the state-of-the-art virtual work energy relative pressure estimator from the 4D-flow measurements, i.e., evaluating formula \eqref{eq:vWERPpdrop} from the synthetically generated data.
\item PINN$+$vWERP: We compute the velocity field from the PINN neural field, thus smoothing and super-sampling the input signal in time. The enhanced velocity field is then plugged-in formula \eqref{eq:vWERPpdrop} to compute the pressure drop.
\item High-order methods: We repeat the test of the three previous methodologies using high-order finite difference schemes for the kinetic power contribution. We anticipate issues with the pure vWERP approach as it relies on a very low time-sampling of the velocity field, which will worsen the estimations as more time information is used for the reconstruction.
\end{enumerate}

\section{Results}\label{sec:results}

In this section, we evaluate the performance of the proposed Physics-Informed Neural Network framework in estimating the velocity, pressure, and viscosity fields within a realistic aortic geometry at five different hematocrit levels. We assess the results both qualitatively and quantitatively by comparing them to high-fidelity FEM simulations, which serve as reference solutions. Our analysis focuses on reconstructed flow structures, mass conservation, and field smoothness, with particular attention given to physiologically relevant metrics such as pressure drops and wall shear stress (WSS). Additionally, we quantify the accuracy of our estimations using statistical measures, including $\ell^2$ relative errors and coefficients of determination ($R^2$), while also identifying variations in reconstruction performance throughout the cardiac cycle.

\subsection{Reconstruction Analysis}

\begin{figure}
    \centering
    \includegraphics[width=0.99\linewidth]{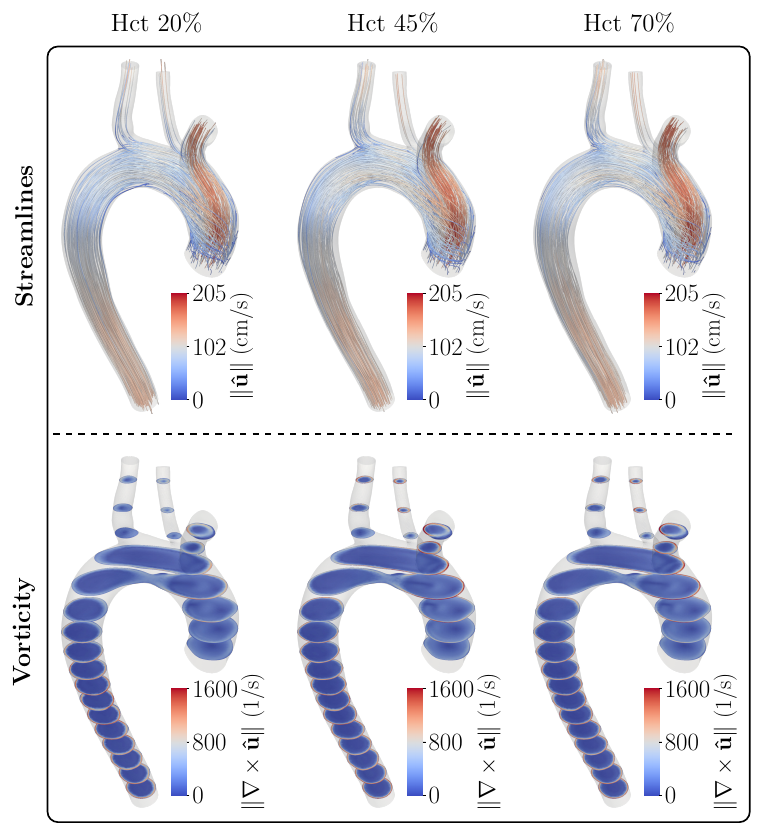}
    \caption{Reconstructed flow fields at peak systole for three hematocrit levels: 20\% (anemic), 45\% (physiological), and 70\% (polycythemic). The first row shows streamlines, revealing well-organized jet-like structures emerging from the inlet and branching into the supra-aortic vessels. The second row presents vorticity magnitude on axial slices, highlighting near-wall shear layers and elevated vorticity in bifurcation regions. The PINN captures characteristic shear-thinning behavior through smooth velocity gradients.}
    \label{fig:stream_vorticidad_syst}
\end{figure}

We begin by analyzing the qualitative behavior of the reconstructed flow fields throughout the cardiac cycle and across different hematocrit levels. Figure \ref{fig:stream_vorticidad_syst} presents the streamlines and vorticity fields at peak systole for three hematocrit levels: low (representing an anemic patient), medium (representing a physiological value), and high (representing a polycythemic patient). The estimated flow fields align well with the measured data, showing well-structured flow patterns. Additionally, the no-slip boundary conditions are accurately represented across all hematocrit levels. The vorticity field remains low in the bulk of the fluid but increases near the boundaries due to the effects of the viscous boundary layer. Notably, the highest vorticity values are observed in the supra-aortic branches. In contrast, the outflow velocity profile resembles a classical shear-thinning pipe solution, with significant vorticity mainly occurring near the walls.

\begin{figure}
    \centering
    \includegraphics[width=0.99\linewidth]{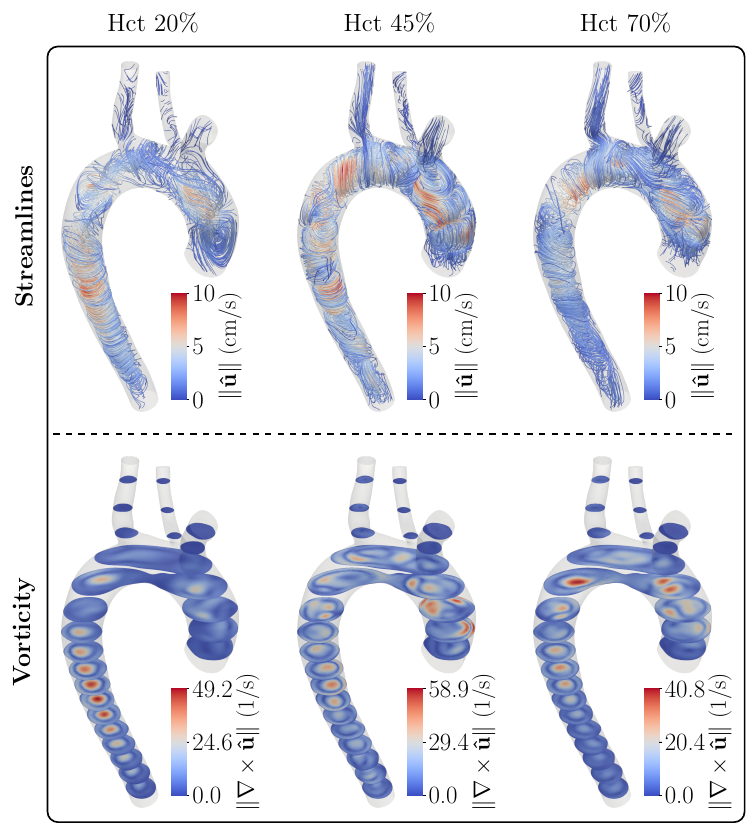}
    \caption{Reconstructed flow fields at mid-diastole for hematocrit levels of 20\%, 45\%, and 70\%. The first row displays streamlines, illustrating complex and fragmented flow patterns with reduced coherence compared to peak systole. The second row shows vorticity magnitude on cross-sectional slices, where multiple small-scale vorticity structures emerge due to the loss of inertial dominance and increased sensitivity to vessel curvature and pressure gradients. Nevertheless, the network preserves key features such as recirculatory regions and near-wall vorticity, which vary in intensity and structure with hematocrit level.}
    \label{fig:stream_vorticidad_dias}
\end{figure}

Figure \ref{fig:stream_vorticidad_dias} illustrates the streamlines and vorticity fields at mid-diastole. When comparing peak systole to mid-diastole, it is clear that blood flow during the diastolic phase is significantly more complex. This complexity arises from a decrease in pressure and the loss of inertial forces. Furthermore, the shape of the aorta contributes to the division of blood flow into smaller vortex structures, which is evident across all three analyzed hematocrit levels. Based on the streamlines and vorticity fields, it can be anticipated that estimating fields related to diastolic conditions will be more challenging, particularly due to the variations in hematocrit levels and their effects on rheological behavior.

\begin{figure}
    \centering
    \includegraphics[width=0.99\linewidth]{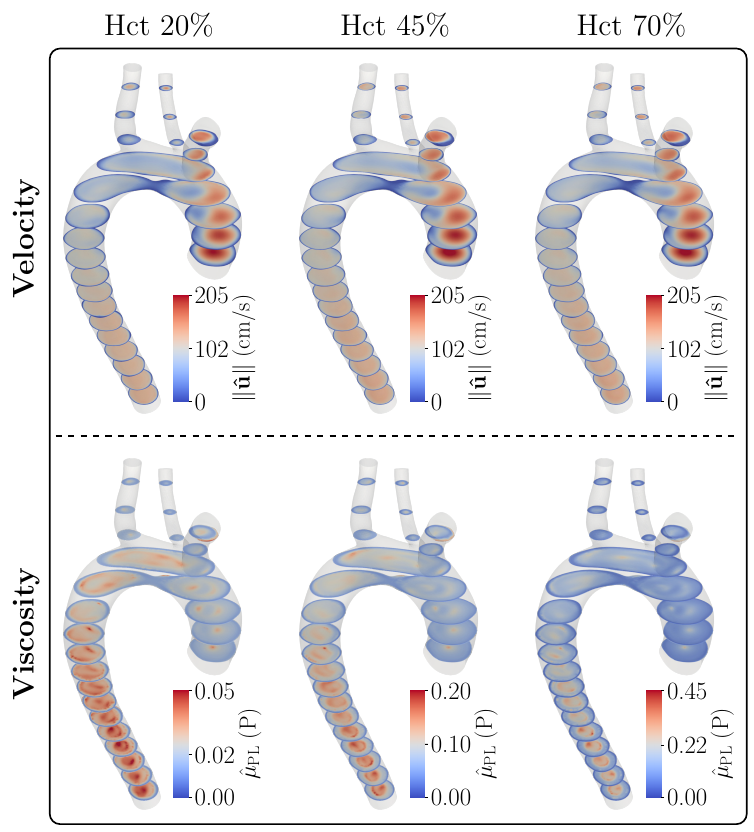}
    \caption{Reconstruction of velocity magnitude, with fully developed flow in the descending aorta. The bottom row displays the corresponding apparent viscosity fields (computed via AD), which reflect shear-thinning behavior. These patterns are consistent with power-law rheology and increase in contrast with hematocrit level.}
    \label{fig:vel_mu}
\end{figure}

Figure \ref{fig:vel_mu} illustrates the estimated velocity magnitude and apparent viscosity fields for various segmentations along the artery. This figure highlights the complexity of the inlet region from a dynamic perspective, which is influenced by bifurcations and changes in flow direction. The velocity magnitude fields reveal an asymmetrical mixed region at the inlet, while a symmetric and fully developed flow condition is observed in the descending aortic branch. The estimated viscosity fields (depicted at the bottom of Figure \ref{fig:vel_mu}) correspond with the velocity fields, indicating shear-thinning regions where the velocity gradient is high and areas of high viscosity where the velocity gradient is low, such as in the central region of the descending aortic branch. This behavior is consistent with the expected viscosity field associated with straight pipe flow.

\begin{figure}[!ht]
    \centering
    \begin{tabular}{@{}c@{}}
        \subfigure[]{\includegraphics[width=0.3\linewidth]{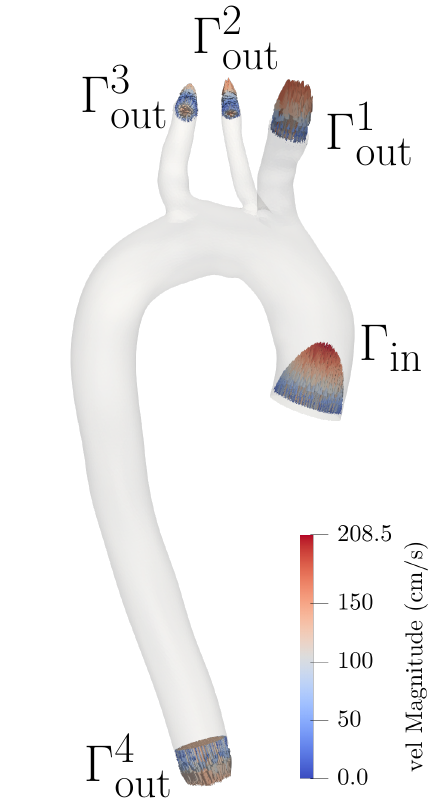}}
    \end{tabular}\qquad
    \begin{tabular}{@{}c@{}}
        \includegraphics[width=0.6\linewidth]{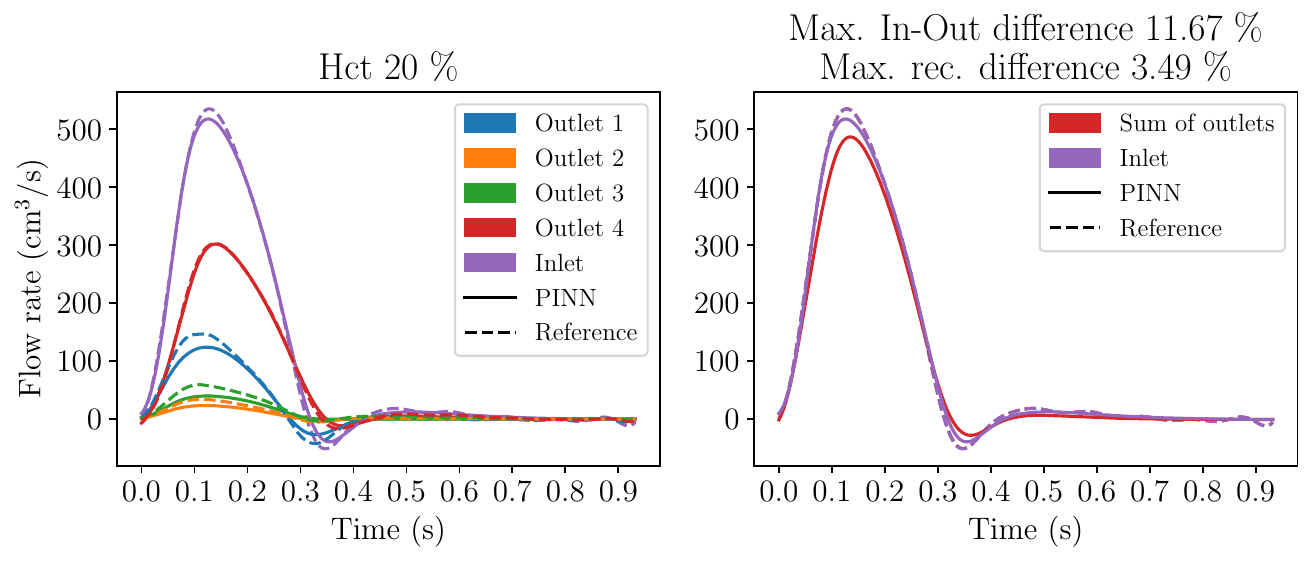} \\
        \includegraphics[width=0.6\linewidth]{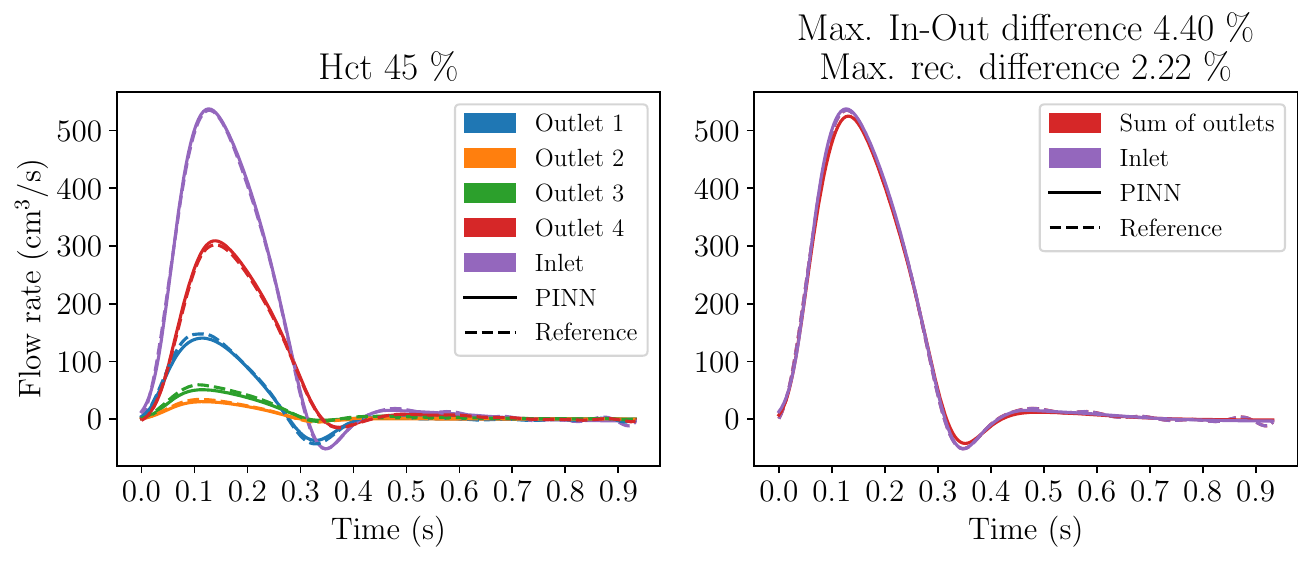} \\
        \subfigure[]{\includegraphics[width=0.6\linewidth]{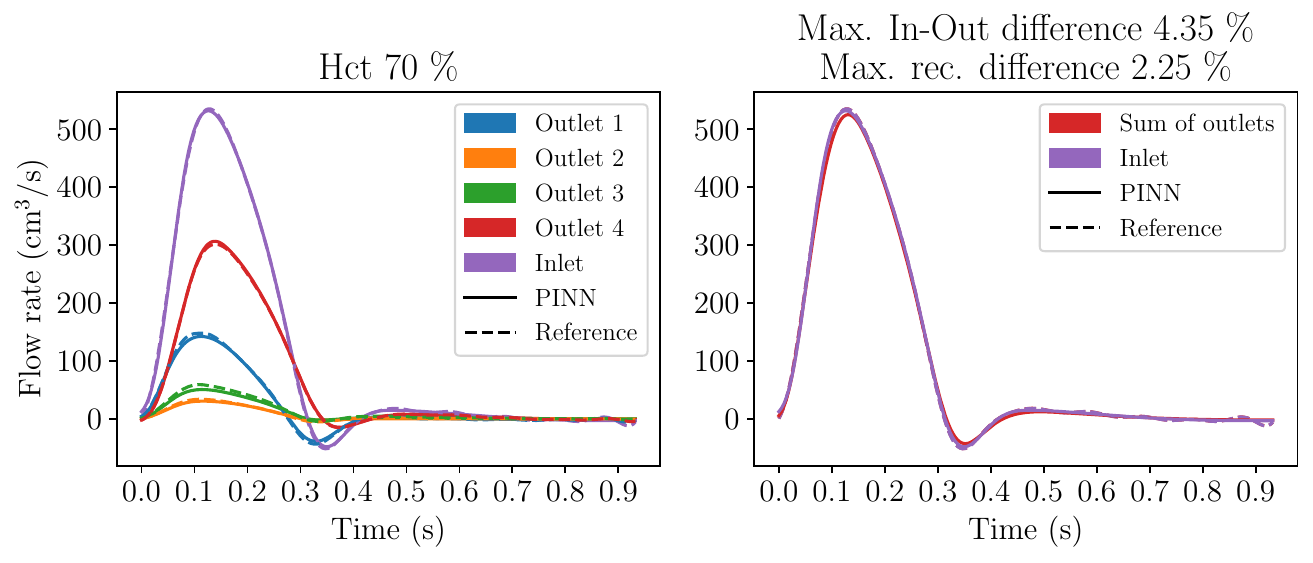}}
    \end{tabular}
    \caption{(a) Position of the cross sections where the flow rates are computed. These are located slightly inwards from the inlet and outlets to avoid the zones with partial volume effects. (b) Flow rates over time in the inlet and outlets. The reported maximum difference between inflow and total outflow is calculated for the PINN reconstruction. This and the maximum reconstruction difference are computed with respect to the maximum reference inflow.}
    \label{fig:flows}
\end{figure}

% P3 flujos: bien
To quantify how the PINN approach describes the bifurcated flows and to evaluate mass flow conservation, Figure \ref{fig:flows} displays the flow rates at the inlet and outlets throughout the cardiac cycle. The highest difference between the inlet and total outlet flows was obtained for the $20\%$ hematocrit, with an error of $11.67\%$ underestimating the outflow. For the $35\%$ hematocrit, the difference was $4.4\%$, and for the $70\%$ hematocrit, $4.35\%$. The error is smaller when comparing the FEM solution to the PINN estimation, with a discrepancy of $3.49\%$ for the low hematocrit case.  %While the estimation is good during systole, it fails to capture oscillations in flow during diastole for all hematocrit levels. 

\begin{figure}[!ht]
    \centering
    \includegraphics[width=0.99\linewidth]{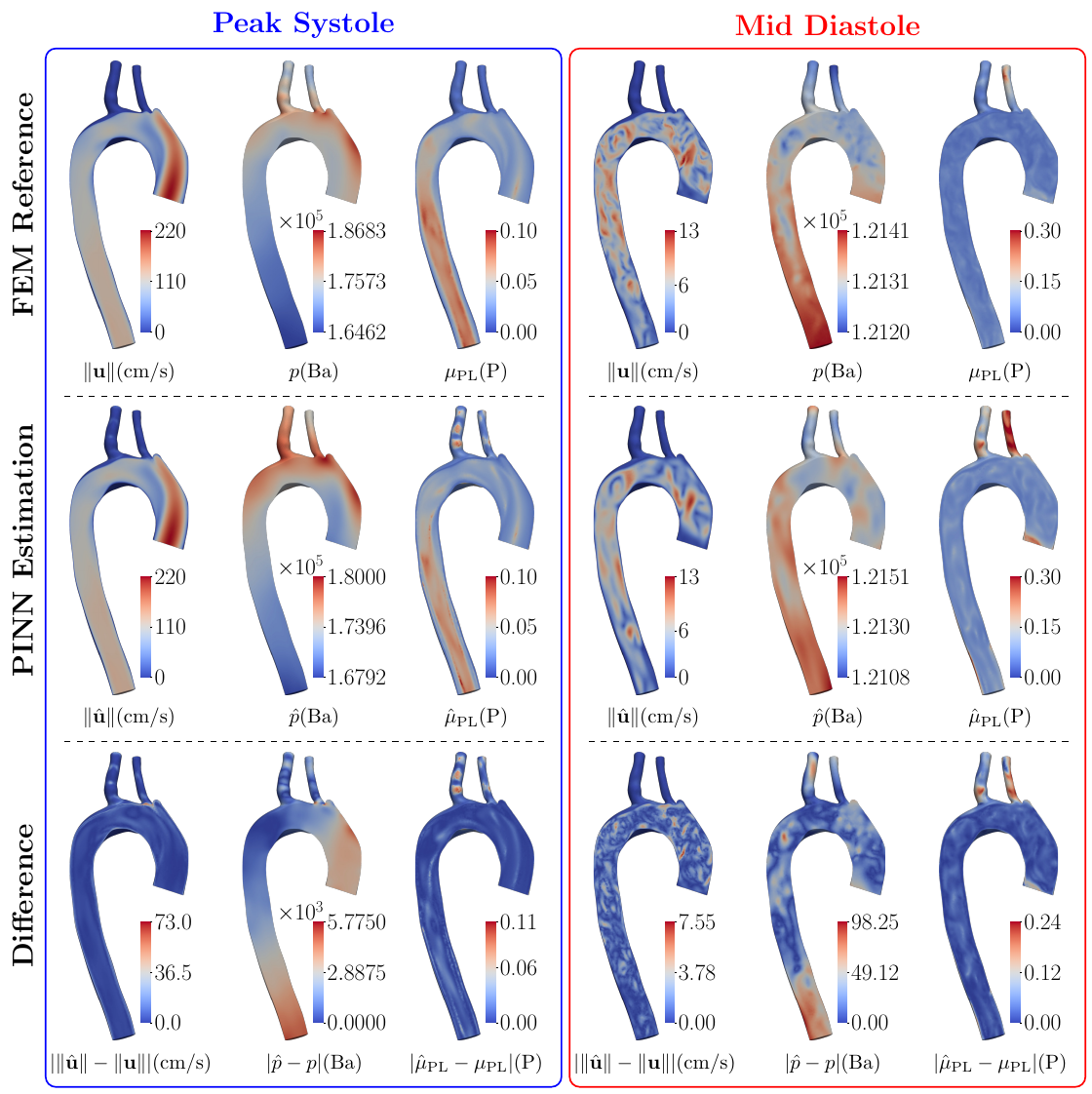}
    \caption{Fields ground truth and PINNs reconstruction for velocity magnitude (in cm/s), pressure (in Ba) and apparent viscosity (in P), for $32.5\%$ hematocrit level, at two time instants $t_1 = 0.13$ s and $t_2=0.75$ s, corresponding to peak systole and mid diastole respectively. The first row shows the reference solutions and the second row the reconstructions. The pressure reconstruction is shifted to have the same mean as the reference at the corresponding instant.}
    \label{fig:velocity_comp}
\end{figure}

Figure \ref{fig:velocity_comp} presents a comparison of the velocity, pressure, and viscosity fields between the estimations of PINN and the FEM solutions. Overall, the reconstruction quality in the full-cut plane is satisfactory. In terms of local maxima and minima values, the PINN approach effectively estimates the velocity and pressure fields. However, the viscosity, which is directly related to the velocity gradients, shows more significant discrepancies. The weak observability of the system for quantities that depend on field derivatives (such as the apparent viscosity) is explained by the fact that the low-fidelity MR measurements provide poor estimates of gradients, so this information is mostly inferred from the physics loss. About flow conditions during peak systole, the differences between PINN and FEM for all analyzed fields are less pronounced compared to those observed during diastole, where more complex fluid dynamics occur.

In diastole, the velocity field lacks detail compared to the reference. We believe that some frequencies present in the MR images associated with the underlying field have not been adequately recovered. This difficulty in finding solutions with high-frequency behaviors, called spectral bias, is a well-known failure mode of PINNs and can be due to a lack of expressivity in the architecture or insufficient training time \cite{Wang2021_ffnet}. The reconstructions for pressure gradients and viscosity are consequently affected.

\subsection{Error Quantification}

\begin{figure}
    \centering
    \includegraphics[width=1\linewidth]{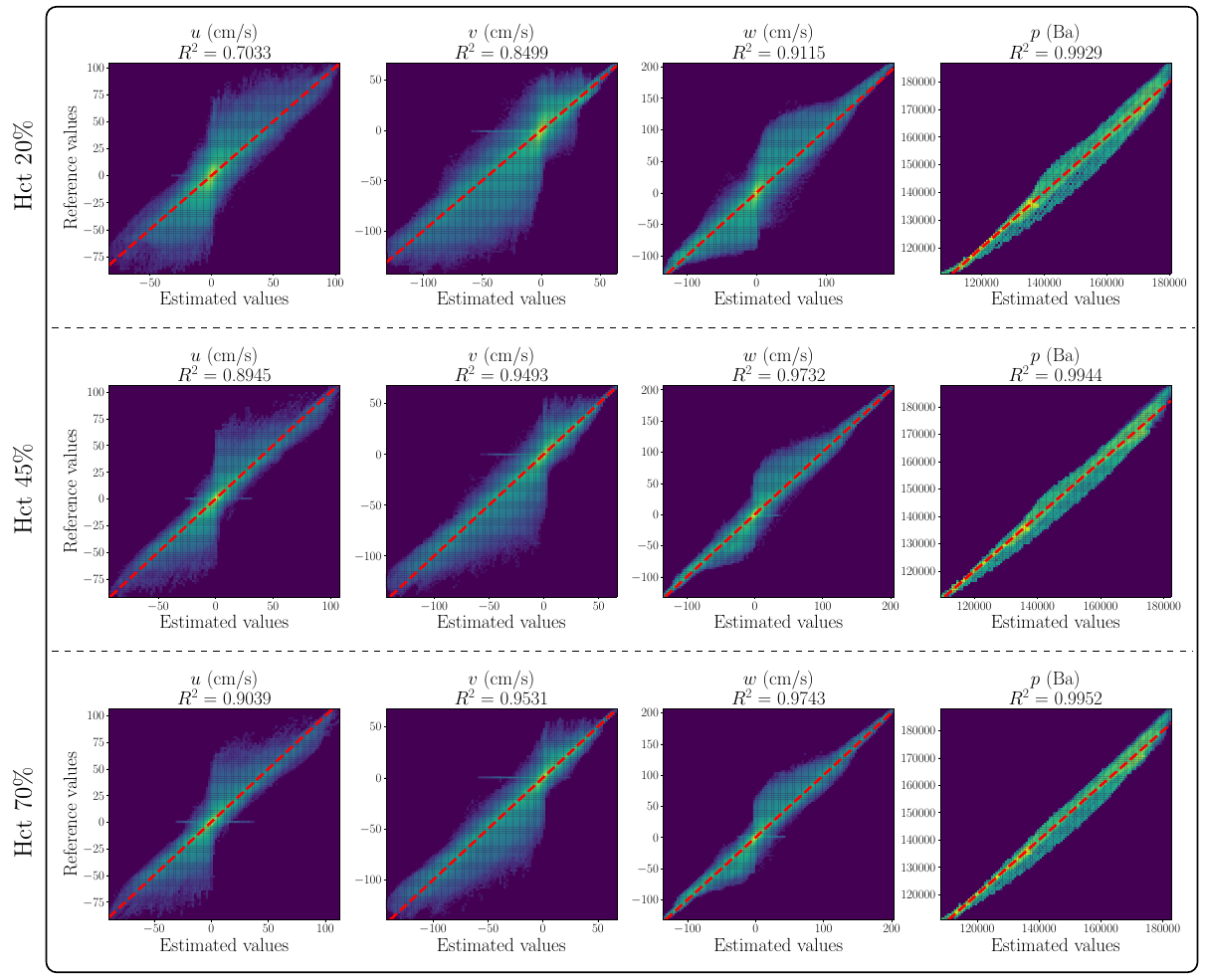}
    \caption{2D Histograms of reference vs. estimated values for the four physical fields and for low, medium and high hematocrit levels. Samples are taken at the nodes of the finite element mesh for $200$ evenly distributed time steps. The heatmap of counts is in a logarithmic scale to aid visualization. Before plotting, the estimated pressure is corrected to have the same average as the reference at each timestep. The coefficient of determination, $R^2$, is computed to quantify the proximity to a perfect estimation (dashed red line).}
    \label{fig:hists}
\end{figure}

This section focuses on quantifying reconstruction errors across all physical fields and time points. To quantify the accuracy of the reconstructed fields, we compare PINN predictions against finite element (FEM) reference solutions over the entire space-time domain. Figure \ref{fig:hists} depicts 2D histograms of reference against estimated values for the fields in the whole space and time domain. The estimation is overall good, with most of the mass concentrated near the perfect-estimation line. However, the model tends to underestimate the velocity components and assign a small velocity to the boundary, resulting in the observed spread in the histogram and a horizontal streak at the zero velocity reference value. Meanwhile, the best estimations are for pressure, without significant spread or any streaks. We quantify the estimation quality using the coefficient of determination $R^2$. Most of the fields shown have an $R^2 \ge 0.9$. Remarkably, the pressure estimations have $R^2 \ge 0.99$. The worst performing examples, the $u$ and $v$ velocity components for the $20\%$ hematocrit level have values of $R^2 =0.7$ and $R^2 =0.85$, respectively.

\begin{figure}[ht]
    \centering
    \includegraphics[width=\linewidth]{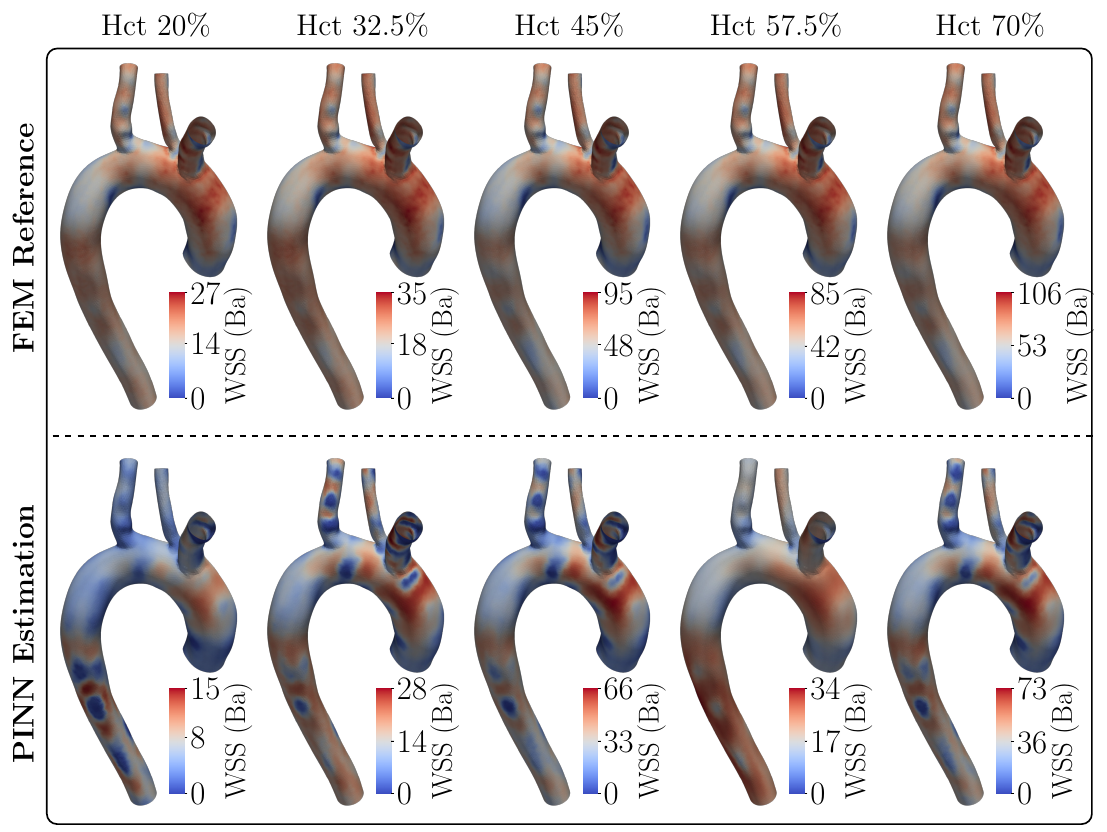}
    \caption{Wall shear stress estimation (in Ba) at the time of maximum value for each hematocrit levels. The first row shows the reference fields, the second row shows the field reconstructions.}
    \label{fig:wss_max}
\end{figure}

As a derived quantity (of interest for the biomedical community), we evaluated wall shear stress over the entire cycle and found the instants with maximum value for all hematocrit levels. By definition, the WSS is the tangential component of the overall boundary viscous force, i.e.:
\begin{equation}
    \mathrm{WSS} = \Vert 2 \mu_{\text{PL}} \left(\mathbf{I} - \mathbf{n}\mathbf{n}^T \right) \left(\dot{\epsilon} \cdot \mathbf{n} \right) \Vert_2,
    \label{eq:wss}
\end{equation}
where $\mathbf{n} : \partial \Omega \rightarrow \mathbb{R}^3$ is a field of unitary normal vectors pointing outward to the working domain. In Figure \ref{fig:wss_max}, we compare the results of the reference FEM solution with those estimated by the PINN approach. While we observe a similar distribution between the two, the PINN reconstructions tend to underestimate the magnitudes when compared to the reference. For instance, at a hematocrit level of $32.5\%$, the reference shows maximum values of approximately $35\,(\mathrm{Ba})$, whereas the reconstructions are around $28\,(\mathrm{Ba})$. Nevertheless, the model successfully reproduces the main features of the reference, including a high-stress region in the ascending aorta. This area is located just upstream of where the inflow jet impacts the wall and disrupts the boundary layer, identifying a critical region of maximum stress.

\begin{figure}[!h]
    \centering
    \begin{minipage}{0.24\textwidth}
        \centering
        \begin{tabularx}{\textwidth}{cc}
            \hline
            Hct (\%) & $\epsilon(p)$ (\%) \\[1ex]
            \hline
            20   & 1.12 \\
            32.5 & 1.13 \\
            45   & 1.00 \\
            57.5 & 1.66 \\
            70   & 0.93 \\
            \hline
        \end{tabularx}
        \captionof{subfigure}{}
    \end{minipage}
    \qquad
    \begin{minipage}{.5\textwidth}
        \centering
        \includegraphics[width=\textwidth]{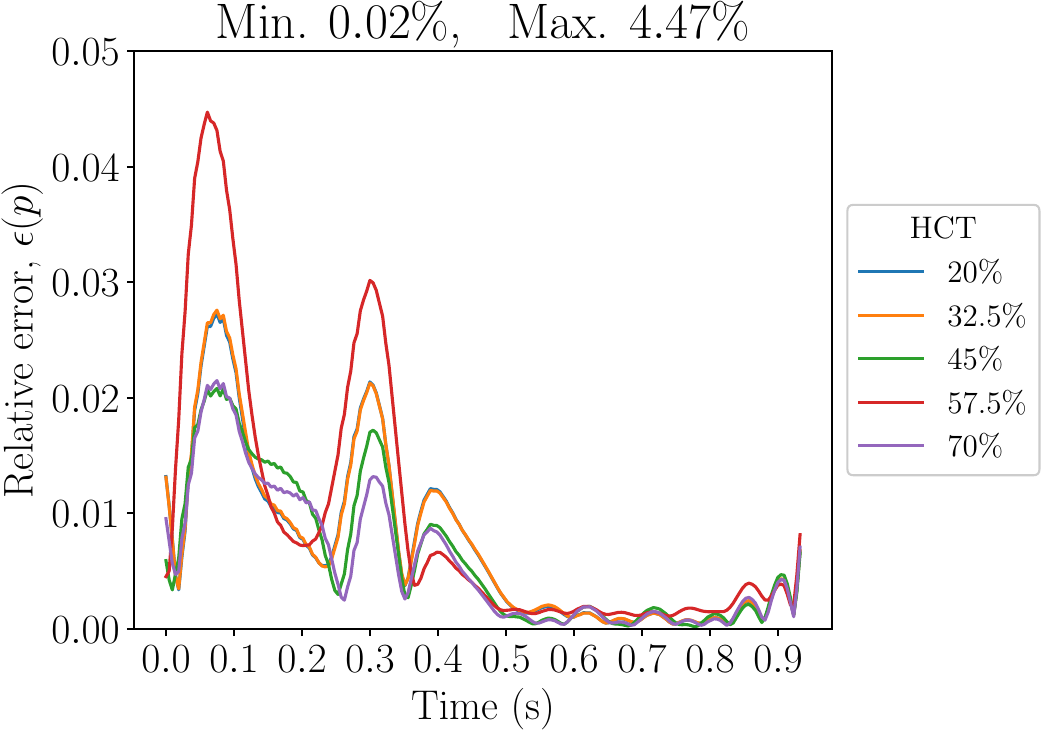}
        \captionof{subfigure}{}
    \end{minipage}
    \caption{(a) Relative error across the space and time domain $\Omega\times[0,T]$, (b) Relative error over the space domain $\Omega$ for $200$ time steps. Errors are computed for the mean corrected pressure at all the hematocrit levels. Minimum and maximum are reported over all time and all hematocrit levels.}
    \label{fig:rel_error}
\end{figure}

We next focus on pressure reconstructions. We compute the relative $\ell^2$ errors as follows:
\begin{equation}\label{eq:rel_error}
    \epsilon(p) = \frac{\norm{p - \hat{p}}_{2}}{\norm{p}_{2}},
\end{equation}
where $p$ is the ground truth FEM solution and $\hat{p}$ stands for the PINN estimation. We compute the relative error over the space-time domain $\Omega\times[0,T]$ and over the space domain $\Omega$ for the $200$ time steps. The estimation is shifted to have the same mean as the reference at each time instant, as is usually done in FEM approaches to enforce zero-mean solutions. Figure \ref{fig:rel_error} shows the reconstruction error for all the hematocrit levels. All cases exhibit similar errors over time, with peaks in systole on the upstroke and downstroke, and errors consistently under $1\%$ during diastole. The highest overall relative error is  $1.66\%$, with a maximum of $4.47\%$ at the systolic upstroke, in the $57.5\%$ hematocrit case. Interestingly, in this case, where the error is greater, the power-law index is the lowest, and the shear-thinning blood behavior is more pronounced. This is a quantitative indicator of the rheological effect and the flow's dependence on the correct estimation of the velocity gradient. Similarly, the low error achieved by the PINN in representing pressure jumps is outstanding, regardless of the hematocrit value computed.

\subsection{Pressure drop estimation}

As described in Section \ref{sec:p_estimation}, we compare four approaches to recover the pressure drop between the input and each of the four outlets. (i) PINN: directly compute the pressure drop from the reconstructed pressure field, (ii) vWERP: compute the vWERP estimation from the MR images, (iii) PINN + vWERP: postprocessing of the PINN reconstructed velocity field applying vWERP, and (iv) repeat the two previous methodologies using a second-order finite difference scheme for the kinetic virtual power contribution. We denote \vWERPfirst{} and \vWERPsecond{} as variants using first- or second-order schemes, respectively. Since the PINN model is smooth in time and inexpensive to query, all estimations are performed at the same time resolution as the FEM reference, that is, 1000 samples over the cardiac cycle, with a time step of $0.937$ ms.

\begin{figure}[!]
    \centering
    \includegraphics[width=1\linewidth]{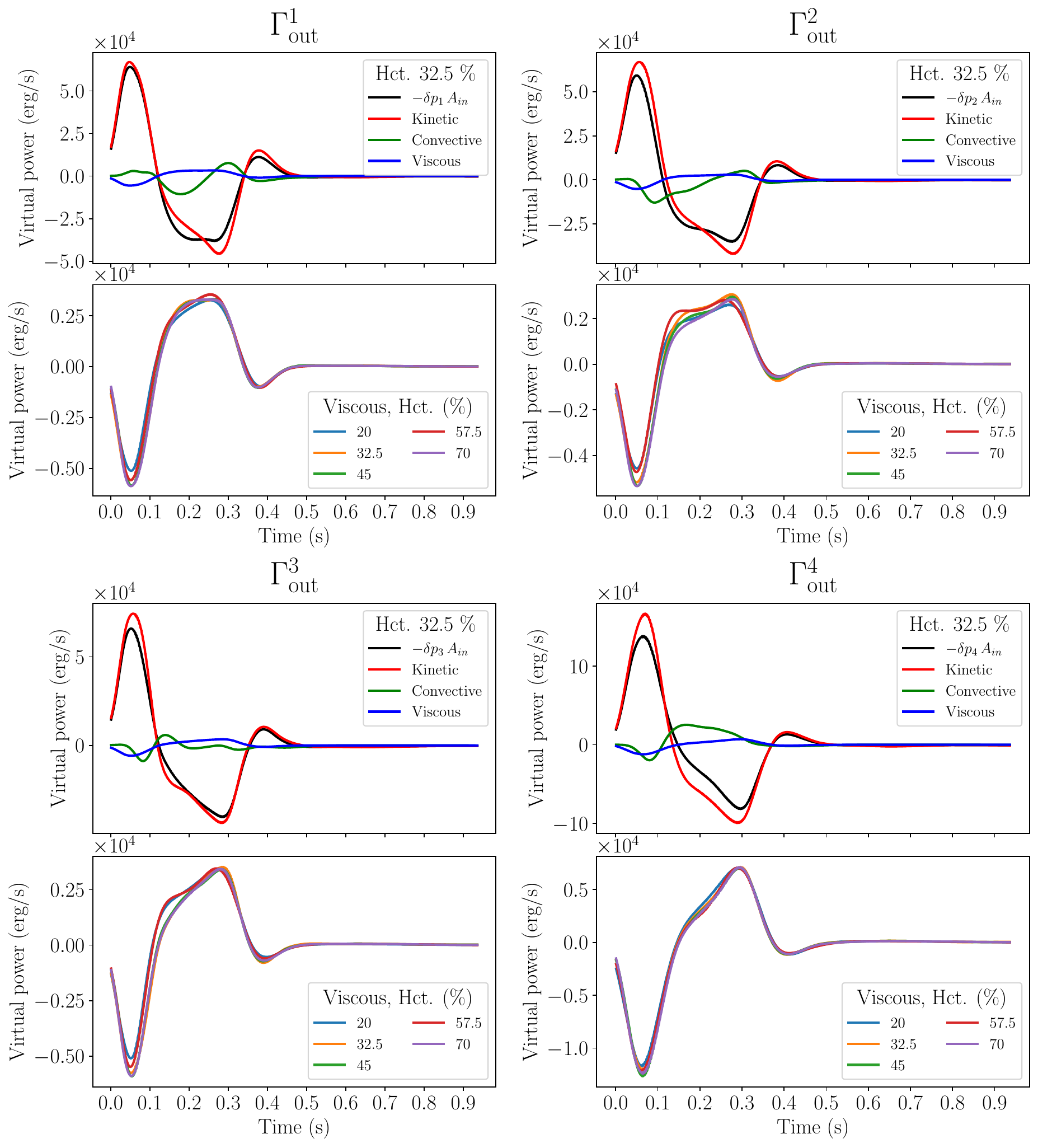}
    \caption{Virtual power computed for each auxiliary test function $\xi^*_k$ ($k=1,\ldots,4$) at fixed hematocrit 32.5 \%. The viscous power is depicted  for all hematocrit levels to highlight its small contribution in the mean outlet pressure in comparison with the kinetic one, which is highly affected by the time-smoothness of the velocity.}
    \label{fig:visc_vWERP}
\end{figure}

Figure \ref{fig:visc_vWERP} presents the virtual power components of the PINN+\vWERPfirst{} approach for a fixed hematocrit level. The kinetic term contributes the most to the pressure drop estimation, while convective and viscous contributions are approximately an order of magnitude smaller, highlighting the importance of velocity time-smoothness to achieve accurate results. The figure also shows the viscous component in all levels of the hematocrit, revealing minimal variation. Therefore, in contrast to other biomarkers, we focus on a single hematocrit level for analysis.

\begin{figure}[!htb]
    \centering
    \includegraphics[width=1\linewidth]{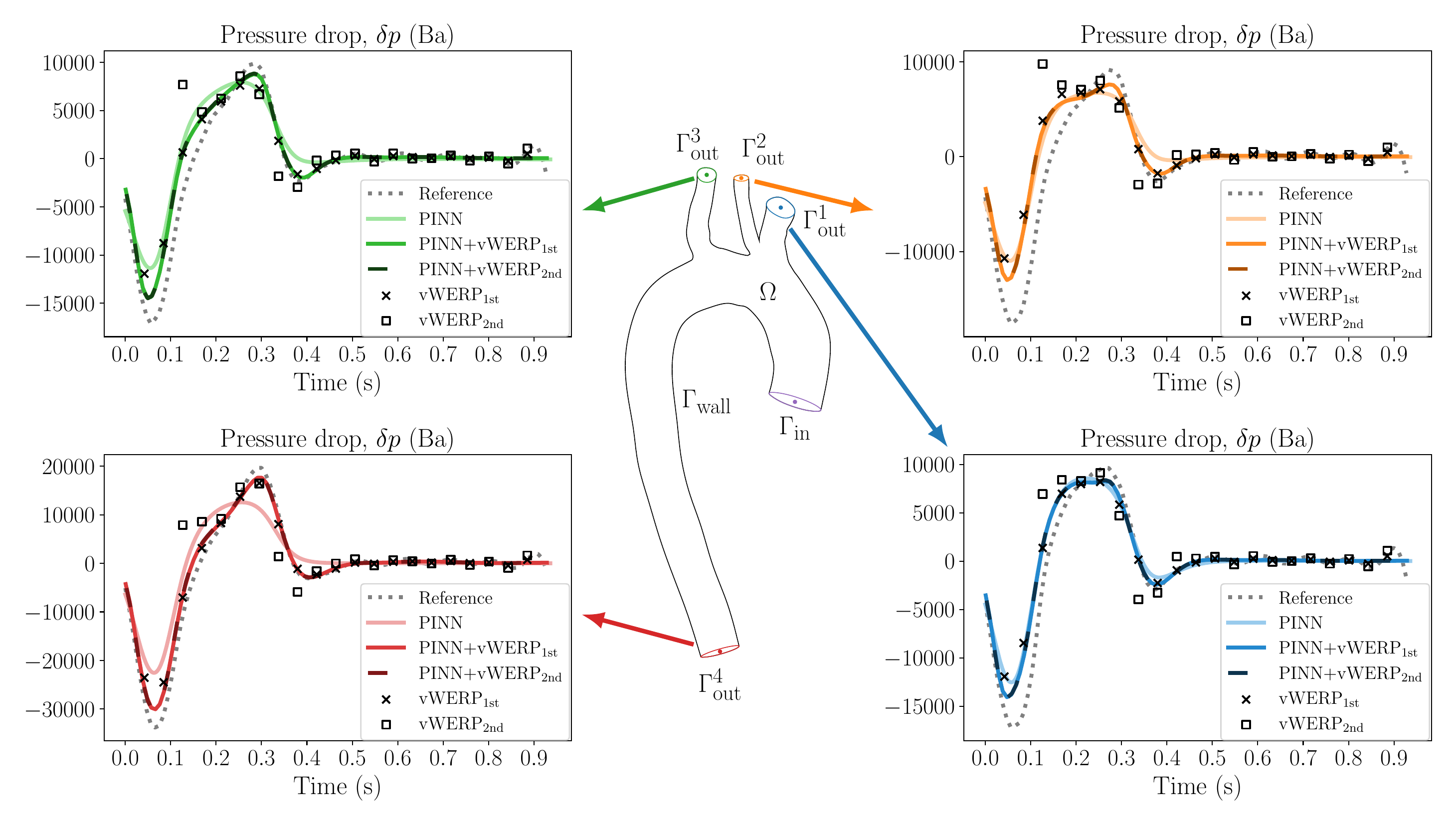}
    \caption{Pressure drop (in Ba) for a fixed hematocrit level of $32.5\%$ between each supra-aortic branch and the ascending aorta. The plots show: the direct result from the reconstructed pressure field, denoted PINN, and the estimations from vWERP applied to the reconstructed velocity field and to the MR images. Both using a first order scheme, denoted PINN+\vWERPfirst{} and \vWERPfirst{}, and using a second order scheme, denoted PINN+\vWERPsecond{} and \vWERPsecond{}.}
    \label{fig:pdrop_32}
\end{figure}

The results of the $32.5\%$ hematocrit level are shown in Figure \ref{fig:pdrop_32}. All approaches accurately predict two peaks, the systolic upstroke and downstroke, with a small temporal offset relative to the reference.
The PINN approach underestimates the peak magnitudes but is otherwise comparable to \vWERPfirst{}. As anticipated, \vWERPsecond{} performs the worst, with larger deviations from the reference during systole.
Additionally, smooth reconstruction by the PINN gives an accurate location of the systolic upstroke and downstroke.
In contrast, vWERP relies entirely on whether an image acquisition coincides with these extrema, which is an unreliable assumption for low time resolutions.
For example, in this study none of the images captures the greatest absolute pressure drop at the systolic upstroke.
Post-processing the velocity reconstruction with vWERP yields the best results, better matching the shape of the reference curve and significantly improving the estimation of peak magnitudes.

To quantify the error in the maximum pressure drop, we define its relative error:
\begin{equation}
    e_{\delta p} = \Bigg\vert\frac{\max\vert\delta \hat p\vert - \max\vert\delta p\vert}{\max\vert\delta p\vert} \Bigg\vert,
\end{equation}
where $\delta\hat p$ and $\delta p$ are the estimated and reference pressure drops, respectively.
Table \ref{tab:diffs_bwn_maxs} reports $e_{\delta p}$ for each method and hematocrit level. In the hematocrit $32.5\%$, direct estimation shows its smallest error ($29.43\%$) for the first outlet. Using PINN+\vWERPfirst{} this error improves to $20.71\%$. Due to the high sampling resolution, PINN+\vWERPsecond{} provides an almost identical result of $20.73\%$. The vWERP post-processing performs best at the fourth outlet, the descending aorta, where the flow is better developed and thus the constant outlet pressure assumption holds best. There, the error is reduced from $34.31\%$ to $11.34\%$ and $10.34\%$, using the PINN+\vWERPfirst{} and PINN+\vWERPsecond{}, respectively.
The largest error, $27.98\%$ and $27.86\%$ occurs for the second outlet, the smallest branch, where the effects of partial volume most significantly impact the MR images and where the pressure distribution is highly uneven. As shown in Table \ref{tab:diffs_bwn_maxs}, these trends persist in all hematocrit levels. The best result is achieved at the level $45\%$ with an error of $8\%$ for the fourth outlet using PINN+\vWERPsecond{}.

\begin{table}[h!]
\centering
\begin{minipage}{0.48\textwidth}
\centering
\caption*{\Large $\Gamma_\mathrm{out}^1$}
\begin{tabular}{cccc}
\toprule
\small Hct & \small PINN & \small +vWERP$_\mathrm{1st}$ & \small +vWERP$_\mathrm{2nd}$ \\
\small (\%) & \small (\%) & \small (\%)                  & \small (\%)                  \\
\midrule
20   & 37.07 & 27.08 & 27.03  \\
32.5 & 29.43 & 20.71 & 20.73  \\
45   & 28.18 & 19.84 & 19.84  \\
57.5 & 44.20 & 23.39 & 23.36  \\
70   & 27.30 & 20.20 & 20.26  \\
\bottomrule
\end{tabular}
\end{minipage}%
\hfill
\begin{minipage}{0.48\textwidth}
\centering
\caption*{\Large $\Gamma_\mathrm{out}^2$}
\begin{tabular}{cccc}
\toprule
\small Hct. & \small PINN & \small +vWERP$_\mathrm{1st}$ & \small +vWERP$_\mathrm{2nd}$ \\
\small (\%) & \small (\%) & \small (\%)                  & \small (\%)                  \\
\midrule
20   & 43.25 & 35.64 & 35.47  \\
32.5 & 39.33 & 27.98 & 27.86  \\
45   & 38.43 & 28.42 & 28.29  \\
57.5 & 64.36 & 36.61 & 36.39  \\
70   & 37.99 & 29.41 & 29.30  \\
\bottomrule
\end{tabular}
\end{minipage}

\vspace{1em}

\begin{minipage}{0.48\textwidth}
\centering
\caption*{\Large $\Gamma_\mathrm{out}^3$}
\begin{tabular}{cccc}
\toprule
\small Hct. & \small PINN & \small +vWERP$_\mathrm{1st}$ & \small +vWERP$_\mathrm{2nd}$ \\
\small (\%) & \small (\%) & \small (\%)                  & \small (\%)                  \\
\midrule
20   & 45.30 & 26.21 & 25.99  \\
32.5 & 35.15 & 17.25 & 17.12  \\
45   & 35.57 & 17.52 & 17.37  \\
57.5 & 62.14 & 23.65 & 23.44  \\
70   & 36.20 & 19.21 & 19.08  \\
\bottomrule
\end{tabular}
\end{minipage}%
\hfill
\begin{minipage}{0.48\textwidth}
\centering
\caption*{\Large $\Gamma_\mathrm{out}^4$}
\begin{tabular}{cccc}
\toprule
\small Hct. & \small PINN & \small +vWERP$_\mathrm{1st}$ & \small +vWERP$_\mathrm{2nd}$ \\
\small (\%) & \small (\%) & \small (\%)                  & \small (\%)                  \\
\midrule
20   & 39.81 & 13.62 & 13.08  \\
32.5 & 34.31 & 11.34 & 10.34  \\
45   & 31.09 &  8.86 &  7.97  \\
57.5 & 69.68 & 12.35 & 11.81  \\
70   & 32.35 & 11.70 & 10.81  \\
\bottomrule
\end{tabular}
\end{minipage}
\caption{Errors in absolute maxima (in \%) for each PINN based method. Each sub table shows the results with respect to each outlet for all hematocrit levels. Each column shows one PINN-based approach.}
\label{tab:diffs_bwn_maxs}
\end{table}

\section{Conclusions}\label{sec:conclusion}

Using a patient-specific aortic model, realistic 4D-flow MRI data in large arteries, and a hematocrit-dependent shear-strain rate relation, we have developed a comprehensive framework for state estimation in blood flows. This problem is significant for the biomedical community due to its potential to non-invasively compute essential quantities for cardiovascular assessment, such as pressure drops and wall shear stresses. 

Our estimates utilize a Physics-Informed Neural Network with several methodological extensions, which, to our knowledge, are being tested for the first time in blood flow applications. The resulting neural network is capable of recognizing significant flow patterns, including fields for velocity, pressure, and viscosity. We observe a satisfactory shear-thinning interpretation, demonstrating that incorporating the governing laws into the optimization problem improves the physical coherence of the fields. Furthermore, we have observed impressive accuracy in estimating pressure drops across both the supra-aortic and descending aortic vessels. We have verified how our approach outperforms the state-of-the-art method vWERP in terms of accuracy and time resolution, as the PINN allows the physical super-sampling of the 4D-flow measurements. In addition, the joint strategy PINN$+$vWERP shows even better results than vWERP or the PINN approach alone. We thus conclude that the best possible approach is to reconstruct first a high-fidelity and physically coherent velocity field from the measurements, and then apply the vWERP estimator to recover the pressure drop with a field that is ensured to have a better Navier-Stokes residual in comparison with the medical data.

Our findings indicate that the neural field can successfully identify critical areas of high wall shear stress using minimal input data. Finally, we note that recovering the diastolic phase is significantly more challenging than recovering the systolic phase. This difficulty comes from the increased complexity of the flow during diastole. It suggests that the network has trouble effectively capturing both types of behavior using a single parameter configuration. This observation aligns with previous work in which optimized models are computed for various time windows throughout the cardiac cycle.

The findings suggest that we continue to develop the tool to better estimate various velocity-related quantities and their derivatives, such as the wall shear stress and the oscillatory shear index. Therefore, we plan to investigate methodological improvements to enhance our results. Implementing domain decomposition could help us address differences in flow dynamics throughout the cardiac cycle and in geometrically distinct areas, such as the inlet, bifurcations, and straight vessels. Additionally, using architectures based on Fourier features may reduce the spectral bias while estimating the diastolic state. Other effective strategies to improve accuracy include employing classical variational methods, using more advanced optimization and training procedures, adjusting variable downstream parameters, and applying mixed multigrid techniques.

% highlights
% - modelo realista (geometria, WKs, y no newtoniano), lo ultimo es poco estudiado
% - nos topamos con el spectral bias en diastole
% - muy poca info de presion fue suficiente para obtener una estimacion decente, no veo razon q tener mas informacion (sin pasarse a algo ficticio) ayude a mejorar la estimacion de presion, esto es ventajoso comparado a garay pero nose como es respecto a otras formas de estimacion de presion.
% - metricas fisicas son importantes, metricas de ML son engañosas
% - un buen estimado de la velocidad no se traduce a buenas estimaciones de campos derivados de esta

% yo tambien, no?. O bien pongo anid magister??
\section*{Acknowledgments} 
M. Sierpe acknowledges the financial support of the Chilean National Agency for Research and Development (ANID) through the program ANID BECAS/MAGISTER NACIONAL/2025-22251449. E. Castillo and F. Galarce acknowledge the financial support of the Chilean National Agency for Research and Development (ANID) through the project Fondecyt No. 1250287. H. Mella acknowledges the financial support provided by ANID through the project Fondecyt No. 11241098. F. Galarce acknowledges DI Vinci PUCV Iniciación 039.731/2025 and the Horizon Europe - 2nd Opportunity OPPTY-MSCA/0125 funding from the Cyprus Research and Innovation Foundation.

\input{main.bbl}
%\bibliographystyle{elsarticle-num}
%\bibliography{literature}

\end{document}

%% file: common.tex
\usepackage[left=1in,right=1in,top=1in,bottom=1in,includefoot,includehead,headheight=13.6pt]{geometry}

\usepackage{xcolor}
\definecolor{darkolivegreen}{rgb}{0.33, 0.42, 0.18}
\definecolor{celestialblue}{rgb}{0.29, 0.59, 0.82}

\usepackage[colorlinks,
            linkcolor=darkolivegreen,
            citecolor=darkolivegreen,
            urlcolor=darkolivegreen]{hyperref}

\hypersetup{pdftitle={CGTJS}, pdfauthor={Felipe Galarce Marin}, 
}

\usepackage{fullpage}
\usepackage{amsmath, amsthm, amssymb}
\usepackage{thmtools, thm-restate}
\usepackage{dsfont}
\usepackage{mathtools}
\mathtoolsset{showonlyrefs=true}
\usepackage{graphicx}   
\usepackage{subfigure}
\usepackage{color}
\usepackage{setspace}
\theoremstyle{definition}

\numberwithin{equation}{section} % Number eqs with section-number

\newcommand{\bP}{\ensuremath{\mathbb{P}}}

\newcommand{\bR}{\ensuremath{\mathbb{R}}}

\def\[{\left[}
\def\]{\right]}
\def\<{\langle}
\def\>{\rangle}
\def\({\left(}
\def\){\right)}
\def\[{\left [}
\def\]{\right]}
\def\({\left(}
\def\){\right)}

\newcommand{\norm}[1]{\Vert #1 \Vert}

\newcommand{\dx}{\ensuremath{\mathrm dx}}
\newcommand{\ds}{\ensuremath{\mathrm ds}}